\documentclass[final,3p,times]{elsarticle}

\usepackage{graphicx}
\usepackage{xcolor}
\usepackage{colortbl}
\usepackage{amssymb,amsmath,amsthm,mathtools,bm}
\usepackage{hyperref}
\usepackage{multirow}
\usepackage{subfig}
\usepackage[section]{placeins}
\usepackage{cleveref}
\usepackage{diagbox}
\usepackage{algorithm}
\usepackage{algpseudocode}
\usepackage{tikz}
\usepackage{booktabs}
\usepackage{float}
\usetikzlibrary{arrows,positioning,fit,backgrounds}
\hypersetup{colorlinks=true,linkcolor=blue,citecolor=blue,urlcolor=blue}

\newcommand{\mm}[1]{{\color{blue}{#1}}}

\AtBeginDocument{%
  \setlength{\abovedisplayskip}{7pt plus 2pt minus 2pt}%
  \setlength{\belowdisplayskip}{7pt plus 2pt minus 2pt}%
  \setlength{\abovedisplayshortskip}{3pt plus 1pt minus 1pt}%
  \setlength{\belowdisplayshortskip}{5pt plus 1pt minus 1pt}%
}

\begin{document}

\begin{frontmatter}

\title{TASE-Stabilized Time-Evolving Natural Gradient Methods for Diffusion-Dominated
PDEs}

\tnotetext[label1]{\mm{Today is \today.} The research of Dongling Wang is supported in part by the National Natural Science Foundation of China {under Grant 12671482} and the Natural Science Foundation of Hunan Province under Grant 2026JJ50363.}

\author[XTU]{Zihao Shi}
\ead{202331510139@smail.xtu.edu.cn}
\author[XTU]{Dongling Wang\corref{cor1}}
\ead{wdymath@xtu.edu.cn}
\cortext[cor1]{Corresponding author.}

\address[XTU]{Hunan Research Center of the Basic Discipline Fundamental Algorithmic Theory and Novel Computational Methods, National Center for Applied Mathematics in Hunan, School of Mathematics and Computational Science, Xiangtan University, Xiangtan, Hunan 411105, China}

\begin{abstract}
The Time-Evolving Natural Gradient (TENG) method evolves neural approximations of time-dependent partial differential equations by projecting time-discrete target states onto the neural-network manifold through local natural-gradient iterations. For stiff diffusion-dominated problems, however, explicit Runge--Kutta methods face severe stability restrictions. High-frequency diffusive components can drive stage targets beyond the locally reachable region of the neural manifold, worsening the conditioning of local least-squares projections and destabilizing parameter updates.

We propose TASE--TENG, which applies the {highly stable explicit operator} $T_p(hW_0)$ to each Runge--Kutta stage increment before projection, where $W_0$ approximates the dominant stiff diffusion operator. To construct this surrogate in a mesh-free setting, we develop a learnable operator whose local interaction coefficients are generated by a shared geometry-dependent neural rule. A structured factorization enforces discrete self-adjointness, dissipativity, and constant preservation. The operator is trained offline and remains fixed during online integration, where it is used only within the TASE stabilization operator.
Numerical experiments on stiff diffusion and reaction--diffusion problems demonstrate improved stability and accuracy at the same time-step size, with reduced errors from high-frequency diffusive modes. The learned rule also allows the stabilization operator to be reconstructed at different spatial point resolutions without retraining.
\end{abstract}

\begin{keyword}
Stiff PDEs \sep TENG method \sep Structure-preserving operator learning \sep TASE stabilization
\end{keyword}

\end{frontmatter}

\section{Introduction}
Neural-network-based methods {have become an active area of research on the numerical solution of differential equations}, including gradient-flow systems, reaction systems, and reaction--diffusion equations. A large class of existing methods follows a global optimization approach, represented by physics-informed neural networks (PINNs) \cite{r1}. In these methods, the governing equation, the initial and boundary conditions, and {additional physical constraints} are included in a single loss functional. Although this formulation is flexible and mesh-free, stiff and multiscale systems often lead to severe optimization difficulties. In particular, stiffness enters directly into the residual loss and typically requires small learning rates to maintain stable optimization. It may also cause competing objectives among different loss terms, ill-conditioned parameter dynamics, and inadequate resolution of fast components \cite{wang2021gradient,krishnapriyan2021failure}. Several recent works {have sought to alleviate these difficulties} by using stiffness-reduction strategies, multiscale network structures, implicit or integral formulations, randomized neural representations, and stiffness-aware sampling mechanisms \cite{Ji2021StiffPINN,Weng2022MPINN,Moya2023DAEPINN,Fabiani2023RanDiffNet,Fabiani2024StabilityRandomProjectionPINN,Zhou2026IPIRNN,Yao2025SolvingMultiscaleDeePODE,kim2025stabilize}.

These developments show that, for stiff dynamics, the main difficulty is not only the approximation power of neural networks. A more important issue is how to represent the stiff physical evolution and how to transfer it stably to the parameter space. This observation motivates time-local neural solvers. Instead of minimizing a global space--time residual, these methods update the neural approximation step by step in time. They naturally follow the causal structure of evolutionary PDEs and can be viewed as neural analogues of classical time-integration methods.

The time-evolving natural gradient (TENG) method \cite{chen2024teng} is a representative time-local framework. At each time step, a time discretization first builds a target state in function space. The neural parameters are then updated by projecting this target onto the neural network manifold through repeated tangent-space approximations. In this way, the PDE evolution is transformed into a sequence of local target-to-manifold projection problems.

{In our previous work \cite{ShiWang2026ConservationTENG}, we combined relaxation of the time-discrete target with projection of the updated neural solution to preserve conservation laws within the TENG framework. The present work addresses diffusion-induced stiffness by stabilizing stage targets and learning a reusable dissipative surrogate for the TASE operator.}

However, for stiff reaction--diffusion equations, the target state generated by an explicit Runge--Kutta method may contain high-frequency diffusive components. These components are often hard to approximate by the local tangent space of the current neural manifold. As a result, the least-squares projection can become ill-conditioned, the parameter updates can become unstable, and {the computation may eventually fail}.

This observation shows that, in stiff regimes, improving TENG should start from the target-construction stage. Instead of forcing the neural manifold projection to fit an unstabilized explicit target, it is better to first construct a target state in which the stiff components have been properly stabilized. In this work, we do this by introducing {Runge--Kutta time-accurate and highly stable explicit (TASE) stabilization} \cite{bassenne2021tase,calvo2021rktase,Conte2026GeneralRKTASE} into the target construction of TENG. More precisely, we apply the TASE stabilization operator $T_p(hW_0)$ to the Runge--Kutta increments, where $W_0$ approximates the main stiff diffusion operator. The resulting stiffness-stabilized target state is more stable and is also more suitable for {the subsequent local neural tangent approximation}.

For strongly stiff reaction--diffusion systems, the main limitation of TENG lies in the stability and local representability of the target state. Explicit Runge--Kutta targets may contain high-frequency diffusive components, and these components are difficult to approximate by the tangent space of the current neural manifold. Therefore, even if the neural network itself is expressive enough, the following target-to-manifold projection may still fail because the parameter manifold cannot locally represent the desired solution update well. This shows that stiffness should be controlled at the target-construction level before the neural parameters are updated. To address this issue, we propose a TASE-stabilized TENG framework. The main contributions of this work are as follows:
\begin{itemize}
\item The Runge--Kutta TASE mechanism is built into the target-construction stage, where the stabilization operator $T_p(hW_0)$ is used to suppress the stiff components in the Runge--Kutta increments. The resulting stabilized target $u_{\mathrm{target}}^{\mathrm{TASE}}$ reduces diffusion-induced high-frequency modes and provides {a target that is more readily represented by the local neural tangent approximation}.

\item A key component of the proposed method is the construction of the stabilization operator $W_0$ required by $T_p(hW_0)$. {To supply the discrete stiff operator needed by the TASE modifier in the neural setting, we construct a learnable dissipative graph surrogate $W_0^{\eta^\ast}$.} Its local coefficients are generated by a shared geometry-dependent neural rule, while dissipativity, {self-adjointness with respect to the discrete inner product}, and constant preservation are enforced analytically through a Gram-type factorization.

After offline training, $W_0^{\eta^\ast}$ is frozen and used only in the TASE modifier. The stabilized Runge--Kutta stage targets are then projected onto the neural manifold in a stage-consistent manner using randomized Stiefel right sketches to improve the stability and robustness of the parameter updates.

\end{itemize}

The rest of this paper is organized as follows. Section~\ref{sec:TimeLocal} reviews the time-local neural approximation framework for evolutionary PDEs and introduces the basic formulation of TENG. Section~\ref{sec:Stiffness}  presents the proposed TASE--TENG method, including the TASE-stabilized target construction, the learnable dissipative graph operator for the stiff diffusion part, and the randomized right-sketched projection strategy for robust parameter updates. Section~\ref{sec:Numerical} {reports numerical experiments on a stiff Prothero--Robinson problem, the viscous Burgers equation, a dryland vegetation model, and the two-dimensional Schnakenberg system to assess the stability and accuracy of the proposed method.} Finally, Section~\ref{sec:Conclusions}  gives concluding remarks and discusses possible future work.

\section{Time-Local Neural Approximation for Evolutionary PDEs}
\label{sec:TimeLocal}
We consider stiff evolutionary PDEs that admit the splitting
\begin{equation}
    \partial_t u(t,\cdot)= \mathcal{F}(u(t,\cdot)):= \mathcal{L}u(t,\cdot)+\mathcal{N}(u(t,\cdot)),\quad u(0,\cdot)=u_0.
\label{eq:stiff_split}
\end{equation}
We assume that $u(t,\cdot)$ belongs to a function space $\mathcal U$ that incorporates the prescribed boundary conditions and permits pointwise evaluation at the sampling points. The operator $\mathcal L$ denotes the stiff dissipative spatial operator, whereas $\mathcal N$ contains the nonlinear reaction and other lower-order terms.

{For the reaction--diffusion systems considered here, let
$u=(u_1,\ldots,u_m)^\top$ denote the $m$-component state, where each
$u_i$ is a scalar-valued function. The Laplacian acts componentwise,
that is, $\Delta u=(\Delta u_1,\ldots,\Delta u_m)^\top$.
With $D=\operatorname{diag}(D_1,\ldots,D_m)$, the diffusion operator
takes the form
$$
\mathcal L u=D\Delta u=
\begin{pmatrix}
D_1 & & 0\\
& \ddots &\\
0 & & D_m
\end{pmatrix}
\begin{pmatrix}
\Delta u_1\\
\vdots\\
\Delta u_m
\end{pmatrix}=
\begin{pmatrix}
D_1\Delta u_1\\
\vdots\\
D_m\Delta u_m
\end{pmatrix},
$$
where $D_i>0$ is the constant diffusion coefficient of the $i$ th
component. This diagonal form describes diffusion without
cross-diffusion terms.
}
Such diffusion operators strongly damp high-frequency spatial modes and therefore constitute a principal source of stiffness in explicit time discretizations.

For scalar or vector-valued functions $f$ and $g$, we define the standard $L^2(\mathcal X)$ inner product and  norm by
$$
    \langle f,g\rangle_{L^2(\mathcal X)} = \int_{\mathcal X} f(x)^{\top}g(x)\,\mathrm dx, \qquad \|f\|_{L^2(\mathcal X)}^2 = \langle f,f\rangle_{L^2(\mathcal X)}.
$$
For vector-valued states, the Euclidean product $f(x)^{\top}g(x)$ naturally includes summation over all state components.
The solution of equation \eqref{eq:stiff_split} is approximated by a neural ansatz $u(t,x)\approx {\hat u}(\theta(t),x)$, $\theta(t)\in\Theta\subseteq\mathbb{R}^{n_\theta}$. This parametrization induces the nonlinear neural approximation manifold $\mathcal{M}= \left\{ \hat{u}_{\theta}: \theta \in\Theta \right\} \subset \mathcal{U}$, $\hat{u}_{\theta}(\cdot)=\hat{u}(\theta,\cdot)$. At the current state $\hat{u}_{\theta}$, the tangent space is
$$
    T_{\hat{u}_{\theta}}\mathcal{M} = \left\{J_{\theta}\xi: \xi\in\mathbb{R}^{n_\theta} \right\}, \, (J_{\theta}\xi)(x) = \nabla_{\theta}\hat{u}^\top\cdot\xi.
$$

{According to the Dirac--Frenkel variational principle\cite{Lubich2008}, the instantaneous PDE vector field is projected orthogonally onto the tangent space of the neural manifold:}
\begin{equation}
    \dot{\theta}(t) \in \arg\min_{\xi\in\mathbb{R}^{n_\theta}} \left\| J_{\theta}\xi - \mathcal{F}(\hat{u}) \right\|_{L^2(\mathcal{X})}^{2} ,
\label{eq:natural_gradient_projection}
\end{equation}
{The first-order optimality condition is}
$$
    \begin{aligned} \left\langle \nabla_\theta \hat u_{\theta(t)}^\top\dot{\theta}(t) - \mathcal F(\hat u_{\theta(t)}), \nabla_\theta \hat u_{\theta(t)}^\top \xi\ \right\rangle = 0, \qquad \forall \xi\in \mathbb{R}^{n_\theta}.\end{aligned}
$$
{The corresponding normal equation is} $G(\theta(t))\dot\theta(t) = b(\theta(t))$, where
$$
    \begin{aligned}
    G_{ij}(\theta) = \left\langle \partial_{\theta_i}\hat u_\theta, \partial_{\theta_j}\hat u_\theta \right\rangle_{L^2(\mathcal X)},\quad b_i(\theta) = \left\langle \partial_{\theta_i}\hat u_\theta, \mathcal F(\hat u_\theta) \right\rangle_{L^2(\mathcal X)}, \qquad i,j=1,\ldots,n_\theta.
    \end{aligned}
$$

The matrix $G(\theta)$ is the tangent-space Gram matrix of the neural manifold at $\hat u_\theta$. The vector $b(\theta)$ represents the $L^2$-projections of the physical vector field $\mathcal F(\hat u_\theta)$ onto the tangent directions. Therefore, the natural-gradient update can be viewed as a least-squares problem that projects $\mathcal F(\hat u_\theta)$ onto the local tangent space $T_{\hat u_\theta}\mathcal M$. However, in practical neural network parametrizations, overparametrization, symmetries, and redundant parameters may make the tangent vectors nearly linearly dependent \cite{finzi2023stable,berman2023randomized}. As a result, $G(\theta)$ can become severely ill-conditioned or nearly singular.

Such degeneracy may render the induced parameter-space evolution system non-unique or poorly determined. Consequently, the resulting parameter update becomes highly sensitive to optimization errors, numerical perturbations, and small changes in the residual projection problem. This issue has motivated a line of work by Benjamin Peherstorfer and collaborators, who proposed randomized sparse-update strategies and Dirac--Frenkel--Onsager-type principles to improve the robustness of parameter-space evolution. {These methods provide mechanisms for stabilizing the tangent-space projection and escaping degenerate or collapsed tangent configurations \cite{berman2023randomized,raviola2026diracfrenkelonsager,Dong2025RandomizedTimeStepping}.}

Since the tangent-space Gram matrix $G$ is often highly ill-conditioned or nearly singular, the parameter-space dynamical system induced by the Dirac--Frenkel principle may be non-unique. Consequently, the continuous evolution in physical space may not admit a direct representation through a stable and smooth trajectory in parameter space. Rather than directly solving the ill-conditioned parameter differential equation, a more robust strategy is to first construct a target state in the physical function space at each time step using a suitable time-discretization scheme and then project this target back onto the neural-network manifold through a local tangent-space approximation. The Time-Evolving Natural Gradient (TENG) method proposed by Chen et al.~\cite{chen2024teng} provides a time-local realization of this idea. It reformulates the continuous projected dynamics as a sequence of target-construction and manifold-projection problems, thereby avoiding the direct long-time integration of an ill-conditioned parameter-space dynamical system.  At time $t_n$, a temporal discretization first constructs a target state in function space,
\begin{equation}
    u_{\mathrm{target}}= \hat{u}(\theta_n,\cdot)+ \Phi_{\Delta t}(\hat{u}(\theta_n,\cdot)),
\label{eq:teng_target_general}
\end{equation}
{where $\Phi_{\Delta t}$ denotes the increment generated by the time-discretization scheme.} For example, the forward Euler target is $u_{\mathrm{target}}= \hat{u}(\theta_n,\cdot)+ \Delta t\,\mathcal{F}(\hat{u}(\theta_n,\cdot))$. TENG then projects this target state back onto the neural manifold by repeated local tangent-space approximations. Starting from $\theta^{(0)}=\theta_n$, each subiteration solves
\begin{equation}
    \Delta\theta^{(k)} \in \arg\min_{\Delta\theta\in\mathbb{R}^{n_\theta}} \left\| u_{\mathrm{target}} -\hat{u} \left( {\theta^{(k)}},\cdot \right)-J_{\theta^{(k)}}\cdot \Delta\theta \right\|_{L^2(\mathcal{X})}^{2},
\label{eq:teng_subiteration}
\end{equation}
followed by $\theta^{(k+1)} = \theta^{(k)}+ \Delta\theta^{(k)}$. {After $K$ subiterations, the next parameter vector is set to} $\theta_{n+1}=\theta^{(K)}$. TENG formulates the time evolution of the PDE as a sequence of local projection problems, in which the target state in the physical space is projected onto the neural-network parameter manifold. At each time step, the temporal discretization constructs a target state in physical space, and the repeated tangent-space projections convert this target into an update of the neural parameters.

{Before applying TENG, we fit the initial condition by finding an initial parameter vector} $\theta_0$ such that the neural network approximation $\hat u(\theta_0,\cdot)$ is close to the given initial condition $u_0$. {This can be done by solving the supervised optimization problem}
$$
    \theta_0= \arg\min_{\vartheta \in \Theta } \frac12 {\left\|\hat u(\vartheta,\cdot)-u_0\right\|^{2}}.
$$

{The initial parameter vector should yield an accurate approximation of the initial condition and a local parametrization suitable for subsequent time integration.} For further details, we refer to the relevant literature \cite{chen2024teng}.

For stiff reaction--diffusion equations, a target generated by an explicit time discretization may contain high-frequency diffusive components. These components are typically poorly aligned with the current neural tangent space, which may lead to ill-conditioned local projections and unstable parameter updates. This motivates the stabilized target construction introduced in the next section.

\section{Stiffness-Aware Target Construction for TENG}
\label{sec:Stiffness}
\subsection{TASE-stabilized Runge--Kutta target construction}
\label{sec:TASE}

TENG can be interpreted as a time-local neural evolution method. Instead of directly evolving the neural parameters through a global space--time optimization procedure, TENG first constructs a target state in the physical space and then projects this target onto the neural manifold through repeated tangent-space approximations. Therefore, the stability of TENG depends on two coupled aspects: whether the physical target construction is stable and whether the resulting target increment can be accurately represented by the local tangent space of the current neural manifold.

For stiff reaction--diffusion equations, explicit Runge--Kutta discretizations may generate stage increments containing strongly amplified high-frequency diffusive components. Although these components are physically dissipative, they may become numerically dominant in the explicit target construction due to the bounded stability region of explicit methods. Consequently, the constructed target state may move away from the current neural manifold, making the tangent-space least-squares projection ill-conditioned and leading to unstable parameter updates.

{To alleviate the diffusion-induced step-size restriction and suppress high-frequency numerical disturbances, we construct a stabilized target $u_{\mathrm{target}}$ using the RK--TASE formulation.}

{
Before incorporating TASE stabilization into the TENG framework, we recall the TASE approach originally introduced by Bassenne,
Fu, and Mani \cite{bassenne2021tase} for the stable numerical integration of stiff differential equations. The presentation
below follows the general RK--TASE formulation developed in \cite{Conte2026GeneralRKTASE}.}
Consider the stiff evolution problem \eqref{eq:stiff_split} after spatial discretization. The main idea of TASE is to introduce a linear stabilization operator $T_p(hW)$, where $h$ is the time step and $W$ is a matrix-form approximation of the dominant stiff component of the evolution operator.
For reaction--diffusion equations, $W$ is typically chosen to approximate the stiff diffusive operator, such as $D\Delta$ or $\operatorname{diag}(D_1,\ldots,D_m)\Delta$. Instead of modifying the Runge--Kutta coefficients, TASE modifies the Runge--Kutta stage increments by applying a rational stabilization operator to the stage right-hand sides.


To explain how TASE alleviates the explicit time-step restriction induced by diffusive stiffness, we first consider the linear semi-discrete system $\frac{\mathrm{d}u}{\mathrm{d}t} = Wu$, where $W\in\mathbb R^{N\times N}$ denotes the spatially discretized diffusion operator. Let $H\in\mathbb R^{N\times N}$ be a symmetric positive definite matrix representing the discrete spatial inner product. For $v,w\in\mathbb R^N$, we define $\langle v,w\rangle_H := v^\top H w$ and $\|v\|_H^2 := v^\top H v$. We assume that $W$ is self-adjoint and dissipative with respect to this inner product:
$$
    W^\top H = HW, \qquad v^\top HWv \le0, \qquad \forall v\in\mathbb R^N.
$$
Consequently, there exists an $H$-orthonormal eigenbasis $\{\phi_k\}_{k=1}^{N}$ such that
$$
    W\phi_k = -\mu_k\phi_k, \qquad \mu_k\geq0, \qquad \phi_i^{\mathrm T}H\phi_j = \delta_{ij},
$$
{where $\delta_{ij}$ is the Kronecker delta. Expanding the solution as} $u(t) = \sum_{k=1}^{N} a_k(t)\phi_k$, {we find that each modal coefficient satisfies} $a_k'(t) = -\mu_k a_k(t)$, $a_k(t) = \mathrm{e}^{-\mu_k t}a_k(0)$. For Laplace-type diffusion operators, the modal decay rates increase quadratically with the spatial frequency. In a Fourier-mode representation, for example, $\mu_k \asymp |k|^2$. More generally, for standard discretizations of second-order diffusion operators, the largest eigenvalue magnitude typically scales as $\mathcal O(h_x^{-2})$ under spatial refinement.

Let $R_{\mathrm{RK}}(z)$ be the stability function of the underlying explicit Runge--Kutta method, and define its negative-real-axis stability radius by
$$
    C_{\mathrm{RK}} = \sup \left\{ C>0: \left|R_{\mathrm{RK}}(z)\right|\leq1 \ \text{for all }z\in[-C,0] \right\}.
$$
For the $k$-th eigenmode, a standard explicit Runge--Kutta step gives $a_k^{n+1} = R_{\mathrm{RK}}(-h\mu_k)a_k^n$. Therefore, stability of all diffusive modes requires
$$
    h\mu_{\max} \leq C_{\mathrm{RK}}, \qquad \mu_{\max} = \max_{1\leq k\leq N}\mu_k.
$$
As the spatial resolution or diffusion coefficient increases, $\mu_{\max}$ grows, and the admissible explicit time step becomes increasingly restrictive. TASE modifies the vector field through the rational operator $\frac{\mathrm{d}u}{\mathrm{d}t} = T_p(hW)Wu$, where
\begin{equation}
    T_p(z) = \frac{\pi_p(z)-z^p}{\pi_p(z)} = 1-\frac{z^p}{\pi_p(z)},
\label{eq:tase_scalar}
\end{equation}
and $\pi_p(z) = z^p - \sigma_1z^{p-1} + \sigma_2z^{p-2} +\cdots+ (-1)^p\sigma_p$. Since $T_p(hW)$ is a rational function of $W$, it commutes with $W$. Thus, the eigenvectors of $W$ remain eigenvectors of the modified operator: $T_p(hW)W\phi_k = -\mu_kT_p(-h\mu_k)\phi_k$. The corresponding RK--TASE modal update is therefore
$$
    a_k^{n+1} = R_{\mathrm{RK}} \left( -h\mu_kT_p(-h\mu_k) \right) a_k^n.
$$
Define the effective stability mapping $q_p(\xi) = \xi T_p(\xi)$, $\xi\leq0$. Using
$\pi_p(\xi) = \xi^p \left( 1-\sigma_1\xi^{-1} +\mathcal{O}(\xi^{-2}) \right), \, \xi\rightarrow-\infty,$
we obtain
$$
    q_p(\xi)=\xi-\frac{\xi^{p+1}}{\pi_p(\xi)}=-\sigma_1+\mathcal O(\xi^{-1}),\qquad \xi\to-\infty.
$$
Hence, $\lim_{\xi\rightarrow-\infty} q_p(\xi) = -\sigma_1$. The unstabilized variable $-h\mu_k$ becomes unbounded as the modal stiffness increases, whereas the TASE-modified variable $q_p(-h\mu_k)$ approaches a finite limit. TASE therefore acts as a rational spectral compression: it maps the unbounded negative diffusive spectrum into a bounded set before the explicit Runge--Kutta stability function is applied.

The condition $\sigma_1 \in [0,C_{\mathrm{RK}}]$ is necessary for stability in the infinitely stiff limit. It is not, however, sufficient for stability at every finite stiffness level, because $q_p$ may exhibit an intermediate-frequency overshoot. For the self-adjoint dissipative problem considered here, a sufficient negative-real-axis stability condition is
$$
    \pi_p(\xi) \neq 0, \qquad \forall \xi\leq0,\qquad q_p\bigl((-\infty,0]\bigr) \subseteq [-C_{\mathrm{RK}},0].
$$
Under these conditions,
$$
    \left| R_{\mathrm{RK}} \left( q_p(-h\mu_k) \right) \right| \leq1 ,
$$
for every $h>0$ and every $\mu_k\geq0$. Consequently, $\|u^{n+1}\|_H \leq \|u^n\|_H$, and the linear diffusive stability is no longer directly constrained by $h\mu_{\max}\leq C_{\mathrm{RK}}$.

We next apply the RK--TASE construction to a general nonlinear or semilinear system $\frac{du}{dt} = \mathcal F(t,u)$. {Let $(A,b,c)$ denote} the Butcher coefficients of an explicit $s$-stage Runge--Kutta method. The RK--TASE stages are defined by
\begin{equation}
    K_i^{\mathrm{TASE}} = hT_p(hW) \mathcal{F} \left( t_n+c_i h,\, u^n+ \sum_{j=1}^{i-1} a_{ij}K_j^{\mathrm{TASE}} \right), \qquad i=1,\ldots,s,
\label{eq:rktase_stage}
\end{equation}
and the final update is
$$
    u^{n+1} = u^n+ \sum_{i=1}^{s} b_iK_i^{\mathrm{TASE}},
$$
Equivalently, every stage can be written as
$$
    \pi_p(hW)K_i^{\mathrm{TASE}} = h \left[ \pi_p(hW)-(hW)^p \right] \mathcal{F}_i,
$$
where
$$
    \mathcal{F}_i = \mathcal{F} \left( t_n+c_i h,\, u^n+ \sum_{j=1}^{i-1} a_{ij}K_j^{\mathrm{TASE}} \right).
$$
Thus, RK--TASE retains the stage coupling of the underlying explicit Runge--Kutta scheme, while the application of the rational stabilization operator introduces one linear solve per stage in the general formulation. In this sense, RK--TASE is more precisely characterized as a linearly implicit stabilization of an explicit Runge--Kutta method \cite{Conte2026GeneralRKTASE}.

The temporal accuracy of RK--TASE methods has been established in \cite{Conte2026GeneralRKTASE}. In particular, if the underlying explicit Runge--Kutta method has order $p$ and the TASE operator satisfies
\begin{equation}
    T_p(hW) = I+C_p(hW)^p+\mathcal{O}(h^{p+1}) = I+\mathcal{O}(h^p), \qquad h\rightarrow0,
\label{eq:tase_consistency_operator}
\end{equation}
then the corresponding RK--TASE scheme preserves the order $p$ of the underlying Runge--Kutta method.
For the rational operator in \eqref{eq:tase_scalar},
$$
    T_p(z) = 1+ \frac{(-1)^{p+1}}{\sigma_p}z^p + \mathcal{O}(z^{p+1}), \qquad z\rightarrow0,
$$
provided that $\sigma_p\neq0$, and hence the required near-identity property follows directly.

Although the stabilization matrix may in principle be constructed from a time-dependent approximation of the stiff Jacobian, repeatedly updating and factorizing such a matrix can be expensive. Following the RK--TASE strategy for diffusion-dominated problems \cite{Conte2026GeneralRKTASE}, we instead employ a fixed stabilization matrix $W_0$ throughout the time integration, $W_n\equiv W_0, \, n=0,\ldots,N_t-1$. For a fixed spatial discretization, a bounded $W_0$ independent of $h$ retains the expansion $T_p(hW_0)=I+\mathcal{O}(h^p)$. In the present setting, $W_0$ is designed to represent the dominant diffusive part responsible for stiffness rather than the complete nonlinear Jacobian.

The stabilized Runge--Kutta stages are then used to construct the physical-space target
$$
    \hat{u}_{\mathrm{target}}^{n+1} = \hat{u}_{\theta^n} + \sum_{i=1}^{s} b_iK_i^{\mathrm{TASE}}.
$$
TENG subsequently projects this target onto the neural approximation manifold,
$$
    \theta^{n+1} \approx \operatorname*{arg\,min}_{\vartheta\in\Theta} \left\| \hat{u}_{\vartheta} - \hat{u}_{\mathrm{target}}^{n+1} \right\|_{L^2(\mathcal{X})}^2.
$$
The RK--TASE consistency result applies to the ideal target constructed in the physical solution space. {Representing this target by the neural ansatz introduces two distinct errors. Let $\mathcal M_{\Theta}=\mathcal M$ denote the neural approximation manifold defined above, and let}
$$
    \theta_{\star}^{n+1} \in \operatorname*{arg\,min}_{\vartheta\in\Theta} \left\| \hat{u}_{\vartheta} - \hat{u}_{\mathrm{target}}^{n+1} \right\|_{L^2(\mathcal{X})}^2,
$$
{be a parameter vector attaining a best approximation of the target on this manifold, assuming that the minimum is attained.} The neural-manifold representation error is defined by
$$
    \varepsilon_{\mathrm{rep}}^{n+1} := \operatorname{dist}_{L^2(\mathcal X)} \left( \hat{u}_{\mathrm{target}}^{n+1}, \mathcal M_{\Theta} \right) = \left\| \hat{u}_{\theta_{\star}^{n+1}} - \hat{u}_{\mathrm{target}}^{n+1} \right\|_{{L^2(\mathcal{X})}}.
$$
This term is intrinsic to the chosen neural architecture and cannot, in general, be eliminated by increasing the number of local projection iterations. The numerical projection error is defined by
$$
    \varepsilon_{\mathrm{proj}}^{n+1} := \left\| \hat{u}_{\theta^{n+1}} - \hat{u}_{\theta_{\star}^{n+1}} \right\|_{{L^2(\mathcal{X})}},
$$
and includes the errors associated with the local tangent linearization, randomized subspace restriction, finite least-squares accuracy, and termination of the TENG iterations. Consequently,
$$
    \left\| \hat{u}_{\theta^{n+1}} - \hat{u}_{\mathrm{target}}^{n+1} \right\|_{{L^2(\mathcal{X})}} \leq \varepsilon_{\mathrm{rep}}^{n+1} + \varepsilon_{\mathrm{proj}}^{n+1}.
$$


\subsection{Construction requirements for \texorpdfstring{$W_0^\eta$}{W0}}\label{subsubsec:w0_requirements}

{The RK--TASE construction in Section~\ref{sec:TASE} requires a reusable discrete approximation of the dominant stiff operator. On a point cloud $X=\{x_i\}_{i=1}^{N_q}$, a sampled stage right-hand side is a vector $r\in\mathbb R^{mN_q}$. Evaluating $T_p(hW_0)r$ requires the action of $W_0$ on sampled vectors and the linear solves associated with the rational modifier. We supply this action through a matrix $W_0\in\mathbb R^{mN_q\times mN_q}$ that remains fixed during the online evolution.}

{Automatic differentiation evaluates spatial derivatives of a given continuous neural function, and hence the PDE right-hand side at the sampling points. However, this evaluation alone does not define a fixed diffusion operator on arbitrary vectors of nodal values: in general, sampled values do not uniquely determine the derivatives of the underlying function. Defining such an operator requires an additional spatial discretization or reconstruction. In our construction, automatic differentiation supplies reference operator actions for offline identification and continues to evaluate the PDE right-hand side during online integration.}

{Since the diffusion coefficients are known, it suffices to construct a scalar Laplace surrogate $L_\eta\in\mathbb R^{N_q\times N_q}$, where $\eta\in\mathbb R^{n_\eta}$ denotes its trainable parameters. For a smooth scalar probe $v:\Omega\to\mathbb R$, define $v_X=(v(x_1),\ldots,v(x_{N_q}))^\top$ and $G_v=((\Delta v)(x_1),\ldots,(\Delta v)(x_{N_q}))^\top$. The identification task is to obtain $L_\eta v_X\approx G_v$ on a representative probe space. With nodal values ordered by component, the full $m$-component surrogate is assembled as}
\begin{equation}
    W_0^\eta = D\otimes L_\eta.
\label{eq:complete_W0}
\end{equation}
{The paired data $(v_X,G_v)$ specify the operator-action fitting problem; the reference action $G_v$ is obtained analytically or by automatic differentiation. The state, tangent, and guard probes and the training objective are specified in Section~\ref{subsec:W0_construction}. After offline training, the selected parameter vector $\eta^\ast$ is frozen, and $W_0^{\eta^\ast}$ is used within the TASE modifier.}

{Operator-action accuracy alone is insufficient for this purpose: a small fitting error on the selected probes need not control the spectrum of an unrestricted matrix. Let $H\in\mathbb R^{N_q\times N_q}$ be the symmetric positive definite matrix defining the discrete scalar inner product. We require $L_\eta^\top H=HL_\eta$, $v^\top HL_\eta v\leq0$ for every nodal vector $v$, and $L_\eta\mathbf1=0$. These identities enforce weighted self-adjointness, dissipativity, and constant preservation independently of training accuracy.}

{Classical finite-difference, finite-element, spectral, or meshfree discretizations can also supply a stiff surrogate; TASE itself does not require a graph neural network. We choose a local graph parametrization to construct the matrix directly from point neighborhoods and to exploit sparse local interactions. A shared geometry-dependent neural rule generates the local coefficients, providing a mechanism for reusing the trained rule across point resolutions. In the construction below, this rule determines first-order difference factors, while a prescribed negative Gram factorization enforces the three structural identities. Thus, the graph specifies locality, the neural map identifies coefficients from operator-action data, and the factorization guarantees the discrete structure. A comparison between the exact and learnable \(W_0\) is provided in \cref{app:learned_exact_W0}.}

\subsection{Structure-preserving graph approximation of the stiff operator}

{The construction proceeds from a local graph on $X$ to first-order difference factors and then to a dissipative matrix factorization. Shared local interactions are standard in graph neural networks; their use in physical simulation is exemplified by MeshGraphNets \cite{Pfaff2021MeshGraphNets}. Here the shared neural map generates the coefficients of the difference factors.}

We next construct a transferable discrete surrogate for the dominant diffusion operator. Consider an $m$-component reaction--diffusion system whose stiff linear part is $\mathcal L U = D\Delta U$, $D = \operatorname{diag}(D_1,\ldots,D_m)$, $D_s>0$.
{
The construction is motivated by the variational structure of the
Laplacian. Let $u$ and $v$ be sufficiently smooth scalar functions
satisfying the prescribed periodic or homogeneous Neumann boundary
conditions. Integration by parts gives
$$
\langle u,\Delta v\rangle_{L^2(\Omega)}
=-\langle\nabla u,\nabla v\rangle_{L^2(\Omega;\mathbb R^d)},
\qquad
\langle v,\Delta v\rangle_{L^2(\Omega)}
=-\|\nabla v\|_{L^2(\Omega;\mathbb R^d)}^2\leq0.
$$
This identity expresses the weak action of the Laplacian through
first-order derivatives and motivates a discrete approximation
based on a negative Gram factorization.
}

{We represent the point cloud \(X\) by a local graph \(\mathcal G_X=(X,\mathcal E_X)\) constructed according to a prescribed neighborhood rule. For each edge \((i,j)\in\mathcal E_X\), the nodal difference \(v_j-v_i\) provides local first-order information for a sufficiently smooth function \(v\). Based on these local differences, we represent the dissipative second-order operator through learnable first-order graph factors \cite{eliasof2021pde}. Unlike standard graph Laplacians with predetermined edge weights, the coefficients of these local factors are generated by a shared geometry-dependent neural map. Consequently, the learned operator is determined by local geometric configurations rather than by a particular point--cloud discretization, allowing the same graph rule to be transferred across different point-cloud resolutions.}

This suggests approximating the second-order stiff operator through a discrete first-order factor.
$$
    v(x_j)-v(x_i) = \nabla v(x_i)^\top (x_j-x_i) + \mathcal O(\|x_j-x_i\|^2).
$$
Therefore, appropriately scaled local graph differences provide a discrete representation of first-order spatial variations.
Based on these local differences, we introduce $R_f$ latent first-order factors. For each factor channel $\ell=1,\ldots,R_f$, we define
\begin{equation}
    \left(Q_\ell^\eta v\right)_i = \sum_{j\in\mathcal N(i)} b_{\eta,\ell}(e_{ij}) \frac{v_j-v_i}{\|x_j-x_i\|}, \qquad \ell=1,\ldots,R_f,
\label{eq:neural_first_order_factor}
\end{equation} where $R_f$ denotes the prescribed number of latent first-order factor channels and $e_{ij}$ denotes a fixed local geometric descriptor associated with the edge $(i,j)$. {The set $\mathcal N(i)$ contains the nodes connected to node $i$ by graph edges.}  The graph topology and the nodal-difference structure are prescribed, whereas the coefficients $b_{\eta,\ell}(e_{ij})$ are learned.

{Equation \eqref{eq:neural_first_order_factor} has a message-passing interpretation \cite{Gilmer2017MPNN}: each summand is an edge message formed from a scaled nodal difference, and summation aggregates the messages at node $i$. The coefficients depend on geometry and $\eta$, so $Q_\ell^\eta$ is linear in the nodal values for fixed $X$ and $\eta$.}

\begin{figure}[!htbp]
\centering
\includegraphics[width=0.95\textwidth]
{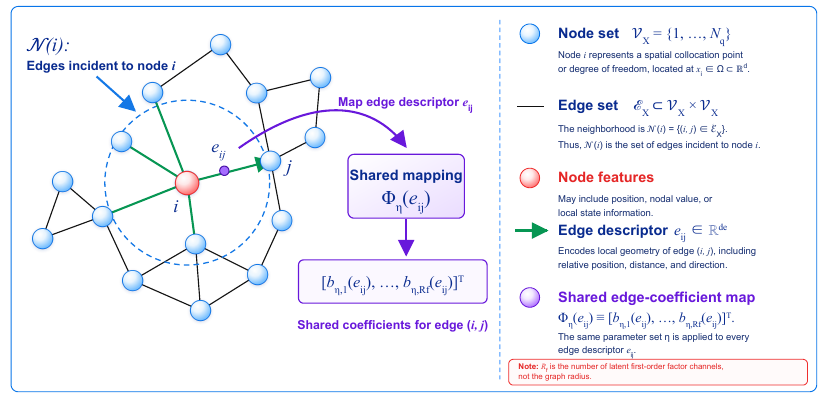}
\caption{{Basic components of the graph structure.}
}
\label{fig:graph_structure}
\end{figure}

The shared coefficient generator is related to the edge-conditioned filter construction of Simonovsky and Komodakis \cite{Simonovsky2017ECC}, in which a neural network maps edge attributes to local filter weights. Here these weights multiply nodal differences, ensuring that each first-order factor annihilates constants. Rather than assigning independent trainable coefficients to individual graph interactions, we use a shared neural map to generate the local coefficients. This avoids a graph-specific parametrization with $\mathcal O(R_f|\mathcal E_X|)$ trainable unknowns and, more importantly, allows the same local rule to be applied to different point-cloud resolutions.

Specifically, let $\Phi_\eta: \mathbb R^{d_e} \rightarrow \mathbb R^{R_f}$, with $\Phi_\eta(e_{ij}) = \bigl(b_{\eta,1}(e_{ij}),\ldots,b_{\eta,R_f}(e_{ij})\bigr)^\top$, where the same parameters $\eta$ are shared by all graph edges.
In the present implementation, $\Phi_\eta$ is represented by a two-hidden-layer neural network. The edge descriptor is first normalized as $\widetilde e_{ij}^{(k)} = \frac{e_{ij}^{(k)}-\mu_k}{\sigma_k}$, $k=1,\ldots,d_e$, and then propagated through
$$
    \begin{aligned} h_{ij}^{(1)}&=\tanh(W_1\widetilde e_{ij}+b_1),\qquad h_{ij}^{(2)}=\tanh(W_2h_{ij}^{(1)}+b_2),\\
     \Phi_\eta(\widetilde e_{ij})&=c_\eta(W_fh_{ij}^{(2)}+b_f)\in\mathbb R^{R_f},\end{aligned}
$$
where $c_\eta>0$ is a trainable global output scale. The $\ell$-th component of $\Phi_\eta$ gives the coefficient $b_{\eta,\ell}(e_{ij})$ used in the $\ell$-th latent first-order factor.

Structure-preserving learning of graph differential operators has been studied through data-driven exterior calculus \cite{Trask2022DDEC} and, more recently, meshfree exterior calculus on point clouds \cite{Shaffer2026MEEC}. The factorization used here is designed to enforce the self-adjointness, dissipativity, and constant preservation required of the TASE surrogate.

The individual factors are stacked as $Q_\eta=\bigl((Q_1^\eta)^\top,\ldots,(Q_{R_f}^\eta)^\top\bigr)^\top$. The scalar discrete Laplace surrogate is then constructed through the Gram-type factorization
\begin{equation}
     L_\eta = -H^{-1}Q_\eta^\top Q_\eta = -H^{-1} \sum_{\ell=1}^{R_f} \left(Q_\ell^\eta\right)^\top Q_\ell^\eta ,
\label{eq:factorized_laplace_surrogate}
\end{equation}

{The graph differences provide a local first-order structure, while training fits the action of the assembled second-order operator $L_\eta$. The channels $Q_\ell^\eta$ are latent factors; the loss does not identify each channel with a particular component of the physical gradient.}
The factorization also preserves the principal structural properties of the Laplace operator. Indeed, $L_\eta^\top H = HL_\eta$, and, for any discrete state $v$,
$$
    v^\top HL_\eta v = -\|Q_\eta v\|_2^2 \leq 0.
$$
{Hence, $L_\eta$ is self-adjoint and dissipative with respect to the $H$-weighted inner product.}  Moreover, because the factors are assembled from nodal differences, $Q_\eta\mathbf 1=0$, $L_\eta\mathbf 1=0$, so that the constant null mode of the periodic or homogeneous Neumann Laplace operator is preserved exactly.

The factorized representation is also consistent with the structure of the target discrete diffusion operator. Let $L_{\eta^\ast}$ denote an ideal discrete approximation of the Laplacian satisfying
$$
    (L_{\eta^\ast})^\top H = HL_{\eta^\ast}, \qquad v^\top HL_{\eta^\ast} v \leq0 .
$$
Then $-HL_{\eta^\ast}$ is symmetric positive semidefinite and therefore admits a Gram factorization $-HL_{\eta^\ast} = (Q_{\eta^\ast})^\top Q_{\eta^\ast}$. Equivalently, $L_{\eta^\ast} = -H^{-1}(Q_{\eta^\ast})^\top Q_{\eta^\ast}$.

Finally, for an $m$-component reaction--diffusion system with $D = \operatorname{diag}(D_1,\ldots,D_m)$, the stiff surrogate used in the TASE operator is defined as
\begin{equation}
    W_0^\eta = \operatorname{diag} \left( D_1L_\eta,\ldots,D_mL_\eta \right).
\label{eq:neural_stiff_surrogate}
\end{equation}
Hence, the construction of $W_0^\eta$ separates three roles: the graph provides the local first-order differential structure, the shared neural map learns how these local interactions are combined, and the Gram factorization recovers the self-adjoint dissipative second-order structure required for TASE stabilization.

\subsubsection{Construction and offline identification of the stiff surrogate}
\label{subsec:W0_construction}

Using the full stiff surrogate $W_0^\eta$ defined in \eqref{eq:complete_W0}, we now establish its discrete structural properties and formulate the offline identification problem.  To represent the discrete inner product for the full $m$-component state space, we introduce the block matrix $\mathbf H_X := I_m\otimes H$, where $I_m\in\mathbb R^{m\times m}$ is the identity matrix. For arbitrary $\zeta,\chi\in\mathbb R^{mN_q}$, define $\langle\zeta,\chi\rangle_{\mathbf H_X} := \zeta^\top\mathbf H_X\chi$. Since $D^\top=D$ and $L_\eta^\top H=HL_\eta$, the complete surrogate is self-adjoint with respect to this block inner product:
$$
    (W_0^\eta)^\top\mathbf H_X = (D\otimes L_\eta^\top)(I_m\otimes H) = D\otimes(L_\eta^\top H) = D\otimes(HL_\eta) = \mathbf H_XW_0^\eta.
$$
To verify dissipativity, let $\zeta=(\zeta_1^\top,\ldots,\zeta_m^\top)^\top\in\mathbb R^{mN_q}$, $\zeta_s\in\mathbb R^{N_q}$, $1\leq s\leq m$. Then,
$$
    \begin{aligned} \left\langle \zeta, W_0^\eta\zeta \right\rangle_{\mathbf H_X} &= \zeta^\top \mathbf H_X W_0^\eta \zeta = \sum_{s=1}^{m} D_s \zeta_s^\top {H}L_\eta \zeta_s = - \sum_{s=1}^{m} D_s \sum_{\ell=1}^{R_f} \left\| Q_\ell^\eta\zeta_s \right\|_2^2 \le0.\end{aligned}
$$
Hence, $\left\langle \zeta, W_0^\eta\zeta \right\rangle_{\mathbf H_X} \le0, \, \forall\zeta\in\mathbb R^{mN_q}$. Therefore, every member of the trainable family $W_0^\eta$ is self-adjoint and dissipative with respect to the $\mathbf H_X$-weighted discrete inner product.

Moreover, since $L_\eta\mathbf 1_{N_q}=0$, the complete surrogate preserves every componentwise constant state: $W_0^\eta \left( C\otimes\mathbf 1_{N_q} \right) = 0$, $\forall C\in\mathbb R^m$. Thus, self-adjointness, dissipativity, and constant preservation are enforced by construction rather than learned from data. The trainable parameters determine only the approximation of the dominant stiff operator within this structurally admissible family.

To further improve the effectiveness of \(W_0^\eta\) when incorporated into the TENG framework, we introduce three complementary families of scalar offline probes, $\mathcal V_{\mathrm{off}} = \mathcal V_{\mathrm{state}} \cup \mathcal V_{\mathrm{tan}} \cup \mathcal V_{\mathrm{guard}}$. {The state probes are the scalar components of the initial neural approximation:} $\mathcal V_{\mathrm{state}} = \left\{ \hat u_{\theta_0,s}: s=1,\ldots,m \right\}$. To represent directions that are locally reachable through parameter variations, we introduce tangent probes
$$
    v_{s,r}^{\mathrm{tan}} = D_\theta \hat u_{\theta_0,s} [\xi_r], \qquad s=1,\ldots,m, \qquad r=1,\ldots,N_{\mathrm{tan}},
$$
where $\xi_r\in\mathbb R^{n_\theta}$ denotes a sampled parameter-space direction.

The guard probes complement the state and tangent probes by constraining the surrogate on functions that are independent of the current neural state. For periodic problems, representative guard probes are paired Fourier modes $\cos(k\cdot x), \, \sin(k\cdot x), \, k\in\mathcal K_{\mathrm{guard}}$. For homogeneous Neumann problems, the guard probes are selected from boundary-compatible cosine modes.

For every $v\in\mathcal V_{\mathrm{off}}$, let $v_X$ and $G_v$ denote the sampled probe and the corresponding reference Laplace action introduced above. The approximation quality of the scalar surrogate is measured by the normalized operator-action error
$$
    \mathcal E_\eta(v) := \frac{ \left\| L_\eta v_X-G_v \right\|_2^2 }{ \|G_v\|_2^2 + \varepsilon_{\mathrm{off}} }, \qquad \varepsilon_{\mathrm{off}}>0.
$$
The normalization prevents probes with large Laplacian magnitudes from dominating the optimization solely because of their amplitudes. The offline training objective is defined by
\begin{equation}
    \begin{aligned} \mathcal J_{\mathrm{off}}(\eta) ={}& \frac{\omega_{\mathrm{s}}} {|\mathcal V_{\mathrm{state}}|} \sum_{v\in\mathcal V_{\mathrm{state}}} \mathcal E_\eta(v) + \frac{\omega_{\mathrm{t}}} {|\mathcal V_{\mathrm{tan}}|} \sum_{v\in\mathcal V_{\mathrm{tan}}} \mathcal E_\eta(v) + \frac{\omega_{\mathrm{g}}} {|\mathcal V_{\mathrm{guard}}|} \sum_{v\in\mathcal V_{\mathrm{guard}}} \mathcal E_\eta(v),\end{aligned}
\label{eq:offline_objective}
\end{equation} where $\omega_{\mathrm{s}}, \omega_{\mathrm{t}}, \omega_{\mathrm{g}} >0$ and $\omega_{\mathrm{s}} + \omega_{\mathrm{t}} + \omega_{\mathrm{g}} = 1$. {The trained surrogate parameters satisfy} $\eta^\ast \in \arg\min_{\eta\in\mathbb R^{n_\eta}} \mathcal J_{\mathrm{off}}(\eta)$. The dependence of the offline objective on $\eta$ follows the differentiable construction
$$
    \Phi_\eta \longrightarrow \left\{ Q_\ell^\eta \right\}_{\ell=1}^{R_f} \longrightarrow L_\eta \longrightarrow \mathcal E_\eta \longrightarrow \mathcal J_{\mathrm{off}}(\eta).
$$
Consequently, gradients of $\mathcal J_{\mathrm{off}}$ with respect to $\eta$ are evaluated by automatic differentiation.

After the offline optimization, the parameter vector $\eta^\ast$ is frozen. We denote the corresponding trained scalar surrogate and full stiff surrogate by
$$
    L_{\eta^\ast} := L_\eta\big|_{\eta=\eta^\ast}, \qquad W_0^{\eta^\ast} := W_0^\eta\big|_{\eta=\eta^\ast} = D\otimes L_{\eta^\ast}.
$$
Since the structural identities above hold for every $\eta\in\mathbb R^{n_\eta}$, they hold in particular at $\eta=\eta^\ast$. Therefore, $\left( W_0^{\eta^\ast} \right)^\top \mathbf H_X = \mathbf H_X W_0^{\eta^\ast}$, and $\left\langle \zeta, W_0^{\eta^\ast}\zeta \right\rangle_{\mathbf H_X} \le0$, $\forall \zeta\in\mathbb R^{mN_q}$. In addition, $W_0^{\eta^\ast} \left( c\otimes\mathbf 1_{N_q} \right) = 0, \, \forall c\in\mathbb R^m$. Consequently, the trained surrogate $W_0^{\eta^\ast}$ preserves the discrete self-adjointness, dissipativity, and constant-preservation properties of the trainable family. In particular, the self-adjointness and dissipativity assumptions required by the linear RK--TASE modal analysis in Section~\ref{sec:TASE} remain valid for the auxiliary linear dynamics generated by $W_0^{\eta^\ast}$.

The structural constraints and the offline approximation objective therefore play complementary roles: the factorized parametrization defines a structurally admissible family $W_0^\eta$, whereas $\mathcal J_{\mathrm{off}}$ selects the member $W_0^{\eta^\ast}$ that best reproduces the dominant stiff operator on the prescribed probe space. The trained operator $W_0^{\eta^\ast}$ is subsequently kept fixed and is used only inside the TASE modifier during the online TASE--TENG evolution; no online retraining of the stiff surrogate is performed.

\subsection{Application of TASE Stabilization to the TENG Framework}
\label{subsec:Application_TASE}

Let $\mathbf U_\theta := \hat u (\theta,X) \in \mathbb R^{mN_q}$ denote the sampled neural state. At time $t_n$, we write $\mathbf U^n := \mathbf U_{\theta_n}$. For the $i$-th stage of an explicit $s$-stage Runge--Kutta method, the TASE-stabilized stage target is defined by
\begin{equation}
    \mathbf Y_i^{\mathrm{TASE}} = \mathbf U^n + \sum_{j=1}^{i-1} a_{ij} \mathbf K_j^{\mathrm{TASE}}, \qquad i=1,\ldots,s.
\label{eq:tase_stage_target}
\end{equation}
The stage target is then represented on the neural approximation manifold. Starting from an initial stage parameter $\theta_{n,i}^{(0)}$, for example $\theta_n$ or the parameter of the preceding stage, the $\ell$-th local TENG iteration evaluates
$$
    J_{n,i}^{(\ell)} := D_\theta\mathbf U_\theta \big|_{\theta=\theta_{n,i}^{(\ell)}} \in \mathbb R^{mN_q\times n_\theta},
$$
and computes a local tangent-space correction from
\begin{equation}
    \Delta\theta_{n,i}^{(\ell)} \in \arg\min_{\Delta\theta\in\mathbb R^{n_\theta}} \left\| \mathbf Y_i^{\mathrm{TASE}} - \mathbf U_{\theta_{n,i}^{(\ell)}} - J_{n,i}^{(\ell)} \Delta\theta \right\|_2^2.
\label{eq:stage_teng_projection}
\end{equation}
The parameter is subsequently updated according to $\theta_{n,i}^{(\ell+1)} = \theta_{n,i}^{(\ell)} + \Delta\theta_{n,i}^{(\ell)}$. In practice, the local correction is computed using the safeguarded randomized right-sketched procedure described in Section~\ref{subsec:robust_projection}.

After $K$ local iterations, we set $\theta_{n,i} := \theta_{n,i}^{(K)}, \, \mathbf U_i^{\mathrm{NN}} := \mathbf U_{\theta_{n,i}}$. To distinguish the neural evaluation of the PDE right-hand side from the learned stabilization operator, we define $F_X(t,\theta) := \left[ \mathcal F \left( t,\hat u_\theta \right) \right](X)$. For the autonomous reaction--diffusion systems considered here, $F_X(\theta) = \left[ \mathcal L \hat u_\theta + \mathcal N(\hat u_\theta) \right](X)$. The spatial derivatives entering $F_X$ are evaluated directly from the continuous neural representation, for example by automatic differentiation. The frozen learned operator $W_0^{\eta^\ast}$ is used only inside the TASE modifier.

The stabilized stage increment is then
\begin{equation}
    \mathbf K_i^{\mathrm{TASE}} = h T_p \left( hW_0^{\eta^\ast} \right) F_X \left( t_n+c_i h, \theta_{n,i} \right), \qquad i=1,\ldots,s.
\label{eq:tase_stage_increment}
\end{equation}
After all stages have been evaluated, the final stabilized target is
\begin{equation}
    \mathbf U_{\mathrm{target}}^{n+1,\mathrm{TASE}} = \mathbf U^n + \sum_{i=1}^{s} b_i \mathbf K_i^{\mathrm{TASE}}.
\label{eq:final_tase_target}
\end{equation}
The final target is then represented on the neural manifold by applying the same local projection procedure, yielding $\mathbf U_{\theta_{n+1}} \approx \mathbf U_{\mathrm{target}}^{n+1,\mathrm{TASE}}$. Thus, the proposed construction is stage-consistent: every Runge--Kutta stage target is first stabilized in the sampled physical space and subsequently represented on the neural manifold before the next stage is evaluated. The TASE modifier controls the stiff components of the physical increment, whereas the local TENG iterations realize the stabilized target in parameter space.

\subsection{Robust Parameter-Manifold Projection}
\label{subsec:robust_projection}

In Section~\ref{subsec:Application_TASE}, the TASE stabilization operator constructs a stiffness-stabilized target $\mathbf U_{\mathrm{target}}^{\mathrm{TASE}}$ in the sampled physical space. The remaining task is to represent this target on the neural approximation manifold through a sequence of local tangent-space corrections.

Two distinct difficulties may arise during this projection. First, although TASE stabilizes the physical-space target before projection, it does not modify the intrinsic conditioning of the neural Jacobian at a fixed parameter state. Parameter redundancy and near-linear dependence among tangent directions may therefore lead to a nearly rank-deficient local parametrization and make the full-space least-squares problem sensitive to numerical perturbations. Second, the target increment may contain components that are only weakly represented by the current tangent space. To obtain a computationally tractable and robust parameter-space realization of the stabilized target, we employ a randomized right-sketched projection.

At the $\ell$-th local TENG iteration, let $\mathbf U^{(\ell)} := \hat u(\theta^{(\ell)},X)$, $J^{(\ell)} := D_\theta \hat u(\theta,X) \big|_{\theta=\theta^{(\ell)}}$, where $J^{(\ell)} \in \mathbb R^{mN_q\times n_\theta}$  {is the Jacobian of the sampled neural state with respect to the parameters.} The corresponding target residual is $\mathbf b^{(\ell)} := \mathbf U_{\mathrm{target}}^{\mathrm{TASE}} - \mathbf U^{(\ell)}$. The full tangent-space correction would then be obtained from
\begin{equation}
    \Delta\theta_{\mathrm{full}}^{(\ell)} \in \arg\min_{\Delta\theta\in\mathbb R^{n_\theta}} \left\| J^{(\ell)}\Delta\theta - \mathbf b^{(\ell)} \right\|_2^2.
\label{eq:full_local_teng_ls}
\end{equation}
For an overparametrized neural representation, solving \eqref{eq:full_local_teng_ls} directly can be computationally expensive and may be sensitive to redundant or nearly dependent tangent directions when $n_\theta$ is large.

{Following the randomized time-stepping approach for nonlinearly parametrized solutions \cite{Dong2025RandomizedTimeStepping}, we restrict each correction to a random low-dimensional subspace of the parameter space.} Let $d_{\mathrm{sk}}\ll n_\theta$ denote the sketch dimension. For each replica $r=1,\ldots,q$, we draw a random Stiefel matrix $\Gamma_r \in \mathbb R^{n_\theta\times d_{\mathrm{sk}}}$, $\Gamma_r^\top\Gamma_r = I_{d_{\mathrm{sk}}}$. Restricting the parameter correction to $\Delta\theta_r^{(\ell)} = \Gamma_r z_r^{(\ell)},$ the reduced coefficient $z_r^{(\ell)}\in\mathbb R^{d_{\mathrm{sk}}}$ is determined by
\begin{equation}
    z_r^{(\ell)} \in \arg\min_{z\in\mathbb R^{d_{\mathrm{sk}}}} \left\| J^{(\ell)}\Gamma_r z - \mathbf b^{(\ell)} \right\|_2^2.
\label{eq:right_sketched_local_ls}
\end{equation}
The corrections obtained from the $q$ independent sketches are then combined as
\begin{equation}
    \Delta\theta^{(\ell)} = \frac{1}{q} \sum_{r=1}^{q} \Gamma_r z_r^{(\ell)}.
\label{eq:averaged_right_sketched_update}
\end{equation}

This construction differs from a coordinate-sparse update. Instead of restricting the correction to selected parameter coordinates, each replica explores a random $d_{\mathrm{sk}}$-dimensional subspace of the full parameter space. Consequently, each reduced least-squares problem involves only $d_{\mathrm{sk}}$ unknowns. Averaging over independent sketches reduces the variability associated with an individual random subspace and provides a more robust correction direction. We emphasize that the averaged correction in \eqref{eq:averaged_right_sketched_update} is used within the repeated local TENG iterations and is not intended to reconstruct exactly the full least-squares solution in \eqref{eq:full_local_teng_ls}.

Because \eqref{eq:right_sketched_local_ls} is obtained from a local linearization of the neural manifold, the resulting correction may still be too large for the tangent approximation to remain accurate. We therefore impose a norm constraint on the averaged correction. Specifically,
$$
    \widetilde{\Delta\theta}^{(\ell)} = \tau_\rho^{(\ell)} \Delta\theta^{(\ell)}, \qquad \tau_\rho^{(\ell)} = \min \left\{ 1, \frac{ \rho }{ \|\Delta\theta^{(\ell)}\|_2 + \varepsilon_{\mathrm{tr}} } \right\},
$$
{where $\rho>0$ specifies the maximum admissible parameter displacement} and $\varepsilon_{\mathrm{tr}}>0$ prevents division by zero. Hence, $\| \widetilde{\Delta\theta}^{(\ell)} \|_2 \le \rho,$ up to the small regularization $\varepsilon_{\mathrm{tr}}$. This safeguard limits the magnitude of the local correction and helps keep the update within the regime in which the tangent-space linearization remains reliable.

The rescaled correction is subsequently evaluated using the nonlinear target-fitting loss
$$
    \mathcal J \left( \theta; \mathbf U_{\mathrm{target}}^{\mathrm{TASE}} \right) := \frac{1}{2} \left\| \hat u(\theta,X) - \mathbf U_{\mathrm{target}}^{\mathrm{TASE}} \right\|_2^2.
$$
Let $\mathcal A = \{\alpha_1,\ldots,\alpha_{N_\alpha}\} \subset(0,1]$ be a prescribed set of damping factors ordered from large to small. We select the first $\alpha^{(\ell)}\in\mathcal A$ satisfying
\begin{equation}
    \begin{aligned} & \mathcal J \left( \theta^{(\ell)} + \alpha^{(\ell)} \widetilde{\Delta\theta}^{(\ell)}; \mathbf U_{\mathrm{target}}^{\mathrm{TASE}} \right) \le (1+\delta_{\mathrm{acc}}) \mathcal J \left( \theta^{(\ell)}; \mathbf U_{\mathrm{target}}^{\mathrm{TASE}} \right),\end{aligned}
\label{eq:line_search_acceptance}
\end{equation} where $\delta_{\mathrm{acc}}\ge0$ specifies the admissible relative increase in the nonlinear target-fitting loss. The accepted parameter update is then $\theta^{(\ell+1)} = \theta^{(\ell)} + \alpha^{(\ell)} \widetilde{\Delta\theta}^{(\ell)}$. If no candidate damping factor satisfies \eqref{eq:line_search_acceptance}, the proposed correction is rejected and {the current parameter vector is retained, with} $\theta^{(\ell+1)} = \theta^{(\ell)}.$

The randomized subspace restriction and the two safeguards play complementary roles. The right-sketched least-squares problems reduce the dimension of the local parameter solve and mitigate the sensitivity associated with redundant parameter directions. The norm-constrained rescaling limits the magnitude of the proposed correction, whereas the nonlinear acceptance test evaluates the correction against the actual target-fitting loss on the neural manifold. Consequently, the TASE stabilization and the randomized projection operate at different levels: TASE stabilizes the physical-space target before projection, while the safeguarded right-sketched iteration provides a robust parameter-space realization of the stabilized target.
\begin{algorithm}[!t]
\caption{Stage-consistent midpoint TASE--TENG}
\label{alg:tase_teng}
\begin{algorithmic}[1]
\Require Time step $h$, final time $T_{\mathrm{End}}$,
number of time steps $N_{t}$ satisfying $N_{t}h=T_{\mathrm{End}}$, neural ansatz $U_\theta$, initial data $U_0$, PDE operator $\mathcal F=\mathcal L+\mathcal N$, collocation set $X$, number $K$ of local TENG iterations, and a precomputed frozen stabilization operator $W_{0}^{{\eta}^{*}}$.
\Ensure Final neural parameter vector $\theta_{N_{t}}$.
\State Fit the initial condition and obtain $\theta_0$ such that $U_{\theta_0}(X)\approx U_0(X)$.
\For{$n=0,1,\ldots,N_{t}-1$}
\State Evaluate the current neural state $U^n=U_{\theta_n}(X)$.
\State Evaluate the first-stage PDE right-hand side $F_1=\mathcal F_X(t_n,\theta_n)$.
\State Compute the TASE-stabilized first-stage increment $K_1 = h\,T_2(hW_{0}^{{\eta}^{*}})F_1$.
\State Construct the midpoint target $Y_{n+\frac12} = U^n+\frac12 K_1$.
\State Project the midpoint target onto the neural manifold $\theta_{n+\frac12} = \operatorname{TENG} \left( Y_{n+\frac12};\theta_n \right)$.
\State Evaluate the midpoint PDE right-hand side $F_2 = \mathcal F_X \left( t_n+\frac{h}{2}, \theta_{n+\frac12} \right)$.
\State Compute the TASE-stabilized second-stage increment $K_2 = h\,T_2(hW_{0}^{{\eta}^{*}})F_2$.
\State Construct the final target $Y_{n+1} = U^n+K_2$.
\State Project the final target onto the neural manifold $\theta_{n+1} = \operatorname{TENG} \left( Y_{n+1};\theta_{n+\frac12} \right)$.
\State Set $t_{n+1}=t_n+h$.
\EndFor
\State \Return $\theta_{N_{t}}$.
\end{algorithmic}
\end{algorithm}

\subsection{Computational Complexity of TASE--TENG}
In this section, we discuss the computational complexity of the TASE--TENG method. We first emphasize that the stabilization operator $W_0^{\eta^\ast}$ used in the present framework is identified and constructed during an offline stage. Once trained, the resulting operator can be stored and repeatedly loaded during the subsequent TENG time evolution, without retraining or updating it at each time step. Therefore, the cost associated with the offline construction of $W_0^{\eta^\ast}$ is a one-time preprocessing cost and does not accumulate with the number of time steps. During the online stage, the additional cost introduced by TASE mainly arises from the stage-wise application of the stabilization operator $T_p(hW_0^{\eta^\ast})$. {The main computational cost of the overall TASE--TENG procedure is associated with} the TENG neural-manifold projection, including Jacobian-related computations and the repeated solution of local least-squares parameter-update problems.

To make the complexity estimate explicit, let $N_q$ denote the number of collocation points and let $m$ be the number of state components. We write $N_x = mN_q$ for the dimension of the sampled physical state. Furthermore, let $n_\theta$ denote the number of trainable parameters in the solution network, $N_t$ the number of time steps, $s$ the number of Runge--Kutta stages, and $K$ the number of local TENG iterations used for each target-to-manifold projection. For the randomized right-sketched projection, $d_{\rm sk}\ll n_\theta$ denotes the sketch dimension and $q$ the number of independent Stiefel sketches.

{The surrogate $W_0^{\eta^\ast}$ is trained only once.} Let $N_{\rm off}=|\mathcal V_{\rm off}|$ denote the total number of offline probes, {$N_{\rm opt}$ the number of optimizer iterations} for training $W_0^{\eta^\ast}$, and $k_g$ the average number of graph neighbors per collocation point. Because the factors $Q_\ell^\eta$ are assembled from local graph interactions, a sparse or matrix-free factor action scales approximately as $\mathcal O(N_q k_g R_d R_f)$. Consequently, the dominant operator-action cost of the offline identification scales approximately as $C_{\rm off} = \mathcal O\!\left( N_{\rm opt}N_{\rm off} N_q k_g R_d R_f\right)$, up to the additional cost of evaluating and differentiating the shared coefficient network $\Phi_\eta$. This cost is incurred only once and is independent of the number of subsequent time steps. Meanwhile, the same trained \(W_{0}^{\eta^\ast}\) can also be used for different time steps.

During the online stage, the learned surrogate $W_0^{\eta^\ast}$ remains frozen and is used only through the stage-wise TASE application $K_i^{\rm TASE} = hT_p(hW_0^{\eta^\ast})F_i .$ Owing to the local factorized construction, one matrix-free application of $W_0^{\eta^\ast}$ scales approximately as $C_W = \mathcal O\!\left( mN_q k_g R_dR_f \right).$

If the TASE rational operator is evaluated using an iterative linear solver with $n_{\rm lin}$ iterations, the corresponding stage-wise cost is approximately $C_{\rm TASE} = \mathcal O(n_{\rm lin} C_W).$

Since $W_0^{\eta^\ast}$ is fixed during time integration, the associated operator structure can be reused at all time steps. The remaining online cost is associated with the TENG projection. Following the complexity analysis of the original TENG method \cite{chen2024teng}, the cost of repeated tangent-space projection can be written as
$$
    C_{\rm TENG} = \mathcal O\!\left( C_{\rm lstsq} N_{\rm it}\frac{T}{\Delta t} \right),
$$ where $C_{\rm lstsq}$ denotes the cost of one least-squares update and $N_{\rm it}$ is the number of local TENG iterations.

In the present stage-consistent formulation, each Runge--Kutta stage requires a target-to-manifold projection. Denoting by $s$ the number of such projections per time step and by $K$ the number of local TENG iterations per projection, the corresponding cost becomes
$$
    C_{\rm TENG}^{\rm online} = \mathcal O\!\left( s K C_{\rm lstsq}^{\rm sk}\frac{T}{\Delta t} \right),
$$ where $C_{\rm lstsq}^{\rm sk}$ denotes the cost of one randomized right-sketched least-squares update used in the present implementation.
Therefore, excluding the one-time offline construction of $W_0^{\eta^\ast}$, the online complexity of TASE--TENG can be summarized as
$$
    C_{\rm online}^{\rm TASE-TENG} = \mathcal O\!\left[ s\left( C_{T_p} + K C_{\rm lstsq}^{\rm sk} \right) \frac{T}{\Delta t} \right].
$$

Combining the one-time offline identification cost with the recurring online evolution cost, the total computational complexity of TASE--TENG can be summarized as $C_{\rm total}^{\rm TASE-TENG} = C_{\rm off} + C_{\rm online}^{\rm TASE-TENG}$, that is,
$$
    C_{\rm total}^{\rm TASE-TENG} = \mathcal O\!\left[ N_{\rm opt}N_{\rm off}N_qk_gR_dR_f + s\left( C_{T_p} + K C_{\rm lstsq}^{\rm sk} \right) \frac{T}{\Delta t} \right],
$$
Here, the first term represents the one-time offline construction of the learned stabilization operator, whereas the second term represents the recurring cost of TASE stabilization and TENG projection during time integration.

\section{Numerical experiments}
\label{sec:Numerical}
In this section, we assess the numerical performance of the proposed TASE--TENG method. All experiments are implemented in Python using JAX \cite{bradbury2018jax} with double-precision arithmetic. Unless otherwise specified, reference solutions are computed using Fourier spectral discretization under periodic boundary conditions. The physical target is constructed using the second-order midpoint RK--TASE scheme. Following \cite{Conte2026GeneralRKTASE}, we use
$$
    T_2(z) = \frac{\pi_2(z)-z^2}{\pi_2(z)} = \frac{5-z}{z^2-z+5}, \qquad \pi_2(z)=z^2-z+5,
$$
corresponding to $\sigma_1=1$ and $\sigma_2=5$. Hence,
$$
    T_2(hW_{0}^{{\eta}^{*}}) = \bigl[(hW_{0}^{{\eta}^{*}})^2-hW_{0}^{{\eta}^{*}}+5I\bigr]^{-1} (5I-hW_{0}^{{\eta}^{*}}).
$$

{Numerical accuracy is measured using the following relative $L^2$ and absolute $L^\infty$ errors:}
$$
    E_{L^2}(t_k) = \left( \frac{ \displaystyle \sum_{i=1}^{n_E} \left| \hat{u}\!\left(\theta_k,x_i^{\mathrm{test}}\right) - u^{\mathrm{ref}}\!\left(t_k,x_i^{\mathrm{test}}\right) \right|^2 }{ \displaystyle \sum_{i=1}^{n_E} \left| u^{\mathrm{ref}}\!\left(t_k,x_i^{\mathrm{test}}\right) \right|^2 } \right)^{1/2},\quad
    E_{L^\infty}(t_k) = \max_{1\le i\le n_E} \left| \hat{u}\!\left(\theta_k,x_i^{\mathrm{test}}\right) - u^{\mathrm{ref}}\!\left(t_k,x_i^{\mathrm{test}}\right) \right|.
$$

To improve computational efficiency, an early-stopping criterion is used in all numerical experiments. Specifically, the TENG iteration is terminated once ${Loss}_{\mathrm{fit}}\le 10^{-13},$ and the computation proceeds to the next stage or time step. Unless otherwise stated, this criterion is used throughout.

\subsection{A stiff Prothero--Robinson test problem}

To assess the behavior of TASE--TENG on a standard stiff dynamical system, we consider the two-component Prothero--Robinson-type problem
\begin{equation}
    \frac{\mathrm d\mathbf y}{\mathrm dt} = A\bigl(\mathbf y-\boldsymbol\phi(t)\bigr) + \boldsymbol\phi'(t), \qquad t\in[0,1],
\label{eq:pr_model}
\end{equation} where $A = \operatorname{diag}(-1,-1000)$, $\boldsymbol\phi(t) = \bigl(\cos t, \sin t\bigr)^\top$. Accordingly, $\boldsymbol\phi'(t) = \bigl(-\sin t, \cos t\bigr)^\top$. The initial condition is prescribed consistently with $\boldsymbol\phi(0)$, $\mathbf y(0) = \boldsymbol\phi(0) = ( 1, 0 )^\top$. Hence, the exact solution is $\mathbf y_{\rm exact}(t) = \boldsymbol\phi(t) = ( \cos t, \sin t )^\top$.

\begin{figure}[!htbp]
\centering
\includegraphics[width=\linewidth]{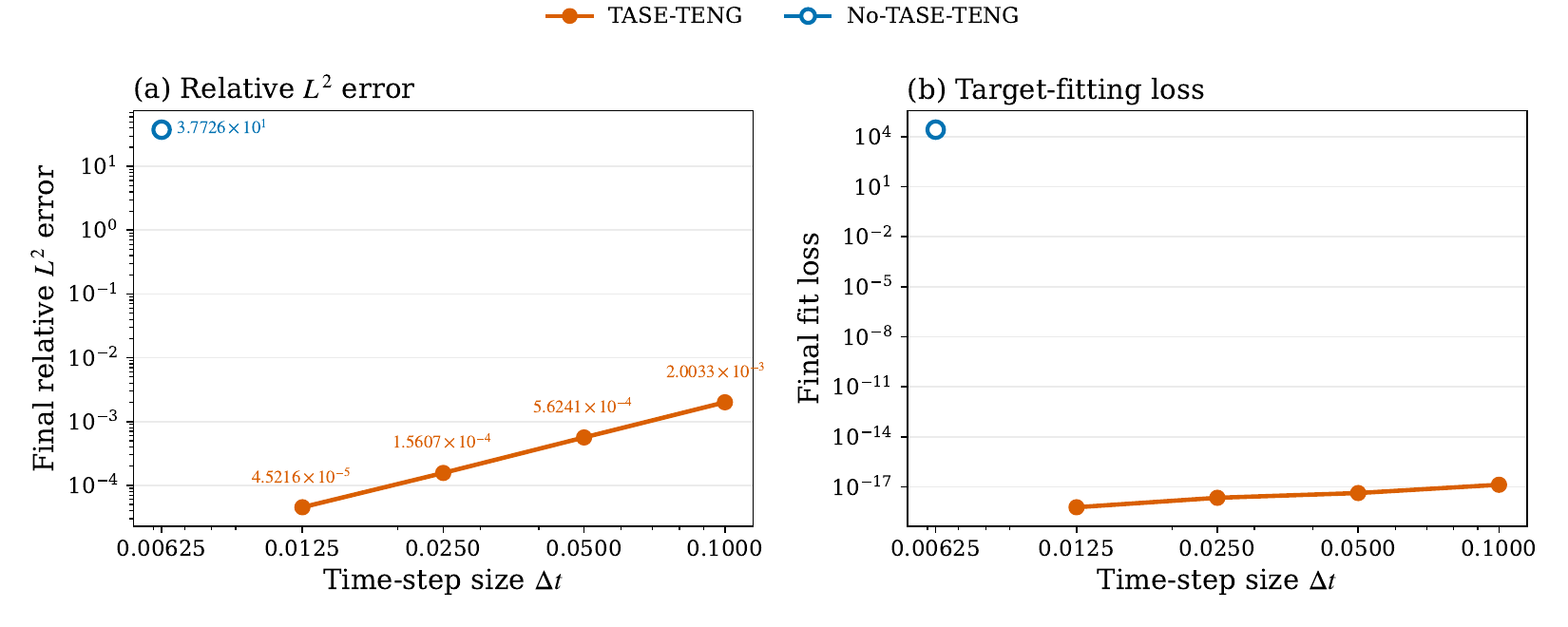}
\caption{Prothero--Robinson problem with $\lambda_{\mathrm{slow}}=-1$ and $\lambda_{\mathrm{stiff}}=-1000$: (a) final-time and trajectory relative $L^2$ errors; (b) final target-fitting loss. Orange curves show TASE--TENG; blue markers show the single No-TASE--TENG run at $\Delta t=0.00625$. The final time is $T=1$, so the number of steps is $1/\Delta t$.}
\label{fig:prothero_robinson_tase}
\end{figure}

{In this stiff ODE test, the two solution components are represented by separate neural networks,} each consisting of two hidden layers of width $8$ with hyperbolic tangent activation functions. Each scalar network contains $97$ trainable parameters. Figure~\ref{fig:prothero_robinson_tase} compares the stabilized TASE--TENG scheme with the corresponding unstabilized midpoint TENG method. The results show that, in this test, {TASE improves the computation by constructing a stable target that can be tracked locally before the TENG projection is performed.} In particular, TASE--TENG remains stable and accurate at $\Delta t=0.1$, whereas the unstabilized counterpart already becomes unstable at the substantially smaller time step $\Delta t=6.25\times10^{-3}$.

\subsection{One-dimensional viscous Burgers equation}

We consider the conservative one-dimensional viscous Burgers equation
\begin{equation}
    \partial_t u = \varepsilon \partial_{xx}u - \partial_x\left(\frac{u^2}{2}\right), \qquad x\in[0,2\pi],
\label{eq:burgers_1d}
\end{equation}
with periodic boundary conditions. The initial condition is prescribed as $u(x,0)=1-\cos (101 x)$ and the viscosity coefficient is chosen as $\varepsilon=0.1$.

The neural solution is represented by a periodic Fourier-feature neural network {to enforce the periodic boundary condition}. More precisely, the input coordinate is embedded through periodic trigonometric features, and the resulting feature vector is passed to a fully connected neural network with three hidden layers and $20$ neurons per layer. The TENG projection is performed on $512$ uniformly distributed collocation points in the interval $[0,2\pi)$.

\begin{table}[!htbp]
\centering
\caption[Error comparison for the viscous Burgers equation]{
Error comparison for the viscous Burgers equation at $T=1.0$. }
\label{tab:burgers_tase_teng_comparison}
\begin{tabular}{lcc}
\toprule
Method & Relative $L^2$ error & $L^\infty$ error \\
\midrule
TASE--TENG & $6.96\times 10^{-5}$ & $1.61\times 10^{-4}$ \\
Standard NG & $1.52\times 10^{-1}$ & $3.76\times 10^{-1}$ \\
Standard NG + random Stiefel & $1.11\times 10^{-1}$ & $2.28\times 10^{-1}$ \\
No-TASE--TENG + random Stiefel & $9.78\times 10^{-2}$ & $2.16\times 10^{-1}$ \\
\bottomrule
\end{tabular}
\end{table}
For comparison, the reference solution is computed using an ETDRK4 Fourier spectral solver with a time step of $\delta t = 10^{-4}$ and $N_{\text{ref}} = 512$ grid points. The TASE--TENG and No-TASE--TENG methods are then tested with a larger time step $\Delta t = 10^{-2}$ up to the final time $T = 1$. This setup is designed to highlight the effect of TASE stabilization on the target construction stage. In the absence of TASE filtering, the explicit midpoint target contains strongly stiff high-frequency diffusive components, which lead to severe oscillations and unstable neural evolution. In contrast, the TASE--TENG solution remains numerically stable and retains a meaningful approximation.

To isolate the effect of TASE stabilization, we compare four method variants against the same reference solution: TASE--TENG, standard NG, standard NG with random Stiefel projection, and No-TASE--TENG with random Stiefel projection. The final-time relative $L^2$ and $L^\infty$ errors are reported in Table~\ref{tab:burgers_tase_teng_comparison}. The results show that TASE--TENG yields substantially smaller errors than the other three unstabilized or only partially stabilized baseline methods. The standard NG method gives a relative $L^2$ error of $1.52\times10^{-1}$. After incorporating random Stiefel projection, this error decreases to $1.11\times10^{-1}$. The No-TASE--TENG variant with random Stiefel projection further reduces the error to $9.78\times10^{-2}$. Nevertheless, these errors remain several orders of magnitude larger than those of TASE--TENG, which attains a relative $L^2$ error of only $6.96\times10^{-5}$ and an $L^\infty$ error of $1.61\times10^{-4}$. This comparison indicates that random Stiefel projection can improve the conditioning of the local parameter update, but it does not by itself provide an accurate stabilized physical target. In contrast, TASE--TENG first constructs a diffusion-stabilized target in the physical solution space and then projects this target onto the neural manifold. Therefore, in this Burgers test, the observed performance gain is primarily attributable to the TASE-based target construction, rather than to the randomized subspace projection alone.

\subsection{Reaction--diffusion vegetation model}

We consider {a reaction--diffusion vegetation model} for plant growth in arid environments, originally introduced in \cite{Eigentler2019Metastability}. The model describes the interaction between two plant species, denoted by $u_1$ and $u_2$, and the available water concentration $w$. Its dynamics are governed by the following three-component reaction--diffusion system:
\begin{equation}
    \left\{ \begin{aligned} \frac{\partial u_{1}}{\partial t} &= \frac{\partial^{2} u_{1}}{\partial x^{2}} +w u_{1}(u_{1}+H u_{2}) -B_{1}u_{1} -S u_{1}u_{2}, \\
     \frac{\partial u_{2}}{\partial t} &= D\frac{\partial^{2} u_{2}}{\partial x^{2}} + F w u_{2}(u_{1}+H u_{2}) - B_{2}u_{2}, \\
     \frac{\partial w}{\partial t} &= d\frac{\partial^{2} w}{\partial x^{2}} + A - w - w(u_{1}+u_{2})(u_{1}+H u_{2}), \end{aligned} \right. \qquad x\in[x_0,x_{\mathrm{end}}],\quad t>0 .
\label{eq:vegetation_model}
\end{equation}

As shown in \cite{Eigentler2019Metastability}, this system admits a metastable coexistence regime for the two species. More precisely, {both species may coexist transiently, whereas one eventually disappears over longer time scales.} Capturing this behavior requires a time-integration method that remains stable and accurate {during long-time simulations}.

In the numerical experiments, we use the parameter values $A=1.5,B_1=0.45,B_2=0.3611,F=0.802,H=0.802,S=0.0002,D=0.802,d=500$. The large diffusion coefficient $d$ in the water equation introduces a strongly dissipative spatial scale, making this problem a suitable benchmark for testing the stability of the proposed TASE--TENG framework. We first consider the time interval $t\in(0,t_{\mathrm{end}}],t_{\mathrm{end}}=4$, on the spatial domain $x\in[-4\pi,4\pi],$ with periodic boundary conditions. The initial data are chosen as $u_1(x,0)=u_2(x,0)=w(x,0)=1+\cos(x)$.

For this three-component system, {we approximate the state components using separate neural networks.} Instead of representing all components by a single shared network, the three state variables $u_1$, $u_2$, and $w$ are approximated by independent scalar subnetworks,
$$
    \hat{U}_{\theta}(x) = \left( \hat{u}_{1,\theta_1}(x), \hat{u}_{2,\theta_2}(x), \hat{w}_{\theta_w}(x) \right), \qquad \theta=(\theta_1,\theta_2,\theta_w).
$$
Each subnetwork is a Fourier-feature MLP with fixed periodic Fourier embeddings followed by fully connected $\tanh$ layers. Specifically, the two vegetation components $u_1$ and $u_2$ use Fourier modes up to $K_u=16$, while the water component $w$ uses modes up to $K_w=4$. All subnetworks have width $20$, depth $3$, Fourier-feature scaling exponent $\beta=0.5$, and an output bias equal to $1$. This separated architecture allows each physical component to have its own parameter space and frequency resolution, while the coupling among $u_1$, $u_2$, and $w$ is imposed through the reaction--diffusion operator during time evolution.

For accuracy assessment, we employ the ETDRK4 method to compute a high-resolution reference solution, with the time step set to $\Delta t = 1\times10^{-4}$ and the number of spatial grid points $N = 512$. Subsequently, the numerical solutions generated by the proposed TASE--TENG method {are compared with ETDRK4}. We set the $R_{f}$ for obtaining $W_{0}$ to 2, and the number of training iterations to 50,000.

\begin{table}[htbp]
\centering
\caption{Relative errors of the vegetation model at \(T=4\) with respect to the reference solution computed using \(\Delta t=10^{-4}\).}
\label{tab:veg_t4_error}
\begin{tabular}{lccccc}
\hline
Method
& \(\Delta t\)
& rel\(L^2(u_1)\)
& rel\(L^2(u_2)\)
& rel\(L^2(w)\)
& rel\(L^2\)(total) \\
\hline

TASE--TENG
& 0.05
& \(3.97\times 10^{-3}\)
& \(1.78\times 10^{-3}\)
& \(4.55\times 10^{-3}\)
& \(3.15\times 10^{-3}\) \\

TASE--TENG
& 0.025
& \(1.58\times 10^{-3}\)
& \(7.20\times 10^{-4}\)
& \(1.87\times 10^{-3}\)
& \(1.26\times 10^{-3}\) \\

TASE--TENG
& 0.0125
& \(5.97\times 10^{-4}\)
& \(2.67\times 10^{-4}\)
& \(7.17\times 10^{-4}\)
& \(4.75\times 10^{-4}\) \\

No--TASE--TENG
& 0.0125
& \(3.33\times 10^{-1}\)
& \(2.11\times 10^{-1}\)
& \(4.62\times 10^{-1}\)
& \(2.84\times 10^{-1}\) \\

TINN \cite{dai2026tinns}
& --
& \(1.09\times 10^{-1}\)
& \(7.92\times 10^{-2}\)
& \(1.97\times 10^{-1}\)
& \(9.76\times 10^{-2}\) \\

\hline
\end{tabular}
\end{table}

To examine the effectiveness of the TASE stabilization operator in the vegetation model under varying time step sizes, we conduct tests using the same random Stiefel embedding design with \(d_{sk}=258\) and \(q=8\). As shown in Table~\ref{tab:veg_t4_error}, the relative \(L^2\) error of TASE--TENG for the water component \(w\) decreases consistently as the time step is refined: from \(4.55\times10^{-3}\) at \(\Delta t=0.05\) to \(1.87\times10^{-3}\) at \(\Delta t=0.025\), and further to \(7.17\times10^{-4}\) at \(\Delta t=0.0125\). In contrast, at the same finest step size \(\Delta t=0.0125\), the method without TASE stabilization (No-TASE--TENG) yields a relative \(L^2\) error as high as \(4.62\times10^{-1}\), which is about two orders of magnitude larger than that of the TASE result. Moreover, the computed profile of No--TASE--TENG at this step size, though still smooth, exhibits severely overestimated amplitudes and a noticeable deviation from the correct physical regime.

To compare TASE--TENG with global space--time training methods, we employ Time-Induced Neural Networks (TINNs) \cite{dai2026tinns} as a baseline. By introducing time-dependent network parameters, TINNs provide adaptive representations at different time instances and are therefore better suited than standard PINNs for time-evolution problems, making them a suitable benchmark for assessing the performance of TASE--TENG in stiff dynamics.

This indicates that, without TASE stabilization, the target increment generated in the explicit Runge--Kutta step remains strongly contaminated by stiff diffusion effects, preventing the subsequent neural--manifold projection from tracking the correct local evolution. Therefore, our experiments confirm that TASE stabilization is essential for constructing reliable local targets in the stiff vegetation problem and significantly enhances the robustness and accuracy of the subsequent TENG projection.
\subsection{ Schnakenberg reaction--diffusion system}
\subsubsection{Periodic boundary conditions}
We consider the two-dimensional Schnakenberg reaction--diffusion system on the periodic domain $\Omega=[0,2\pi)^2$ given by
$$
    \left\{ \begin{aligned} \partial_t u &= D_u\Delta u + \gamma (a-u+u^2v), \\
     \partial_t v&= D_v\Delta v+ \gamma\left(b-u^2v\right).
     \end{aligned} \right.
$$

In the experiments, we use $a=0.1,b=0.9,D_u=1,D_v=100,\gamma=1$. The large diffusion coefficient $D_v$ introduces strong stiffness in the $v$-component, making this problem a suitable benchmark for testing the stability of TASE--TENG. For the neural representation, {we use separate neural networks for the two state components.} Instead of approximating both components by a single shared output network, the two state variables are represented by separate subnetworks:
$$
    \hat{U}_{\theta}(x,y) = \left( \hat{u}_{\theta_u}(x,y), \hat{v}_{\theta_v}(x,y) \right), \qquad \theta=(\theta_u,\theta_v).
$$
Each component is modeled by a fixed Fourier-feature MLP. The Fourier frequencies are controlled separately for the two variables: the $u$-network uses $K_u=4$, while the $v$-network uses $K_v=2$. In the implementation, both subnetworks have width $20$, depth $4$, Fourier-feature scaling exponent $\beta=0.5$, and a small output-layer initialization scale $10^{-3}$. The Schnakenberg system admits the spatially homogeneous equilibrium $u_\ast=a+b, v_\ast=\frac{b}{(a+b)^2}$. For the parameters $a=0.1$ and $b=0.9$ considered here, this gives $(u_\ast,v_\ast)=(1,0.9)$.

The initial condition is chosen as a small periodic perturbation of this equilibrium:
$$
    u_0(x,y) = u_\ast+10^{-2}\cos(x)\cos(y), \qquad v_0(x,y) = v_\ast+10^{-2}\cos(2x)\cos(y).
$$
To quantify the numerical accuracy, we use the root-mean-square error (RMSE) with respect to the ETDRK4 reference solution. For $q\in\{u,v\}$, the componentwise RMSE at the final time is defined as
$$
    \mathrm{RMSE}_{q} = \left[ \frac{1}{N_xN_y} \sum_{i=1}^{N_x} \sum_{j=1}^{N_y} \left| q_{ij}-q_{ij}^{\mathrm{ref}} \right|^2 \right]^{1/2}, \qquad q\in\{u,v\},
$$
where $q_{ij}$ and $q_{ij}^{\mathrm{ref}}$ denote the TENG solution and the reference solution, respectively, evaluated on the same $N_x\times N_y$ grid. For the coupled system, we report the  $\mathrm{RMSE}_{\mathrm{joint}} = \left[ \frac{ \mathrm{RMSE}_{u}^{2} + \mathrm{RMSE}_{v}^{2} }{2} \right]^{1/2}$.

\begin{figure}[!htbp]
\centering
\includegraphics[width=\linewidth]{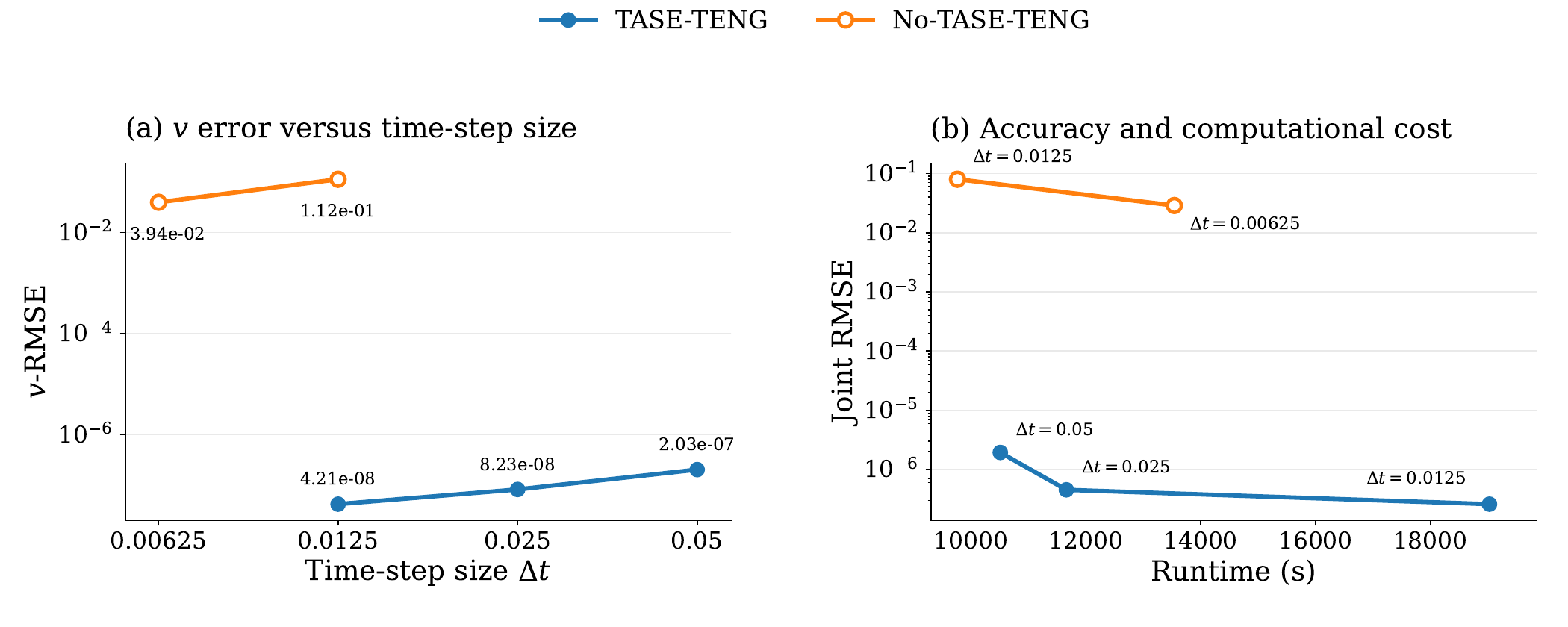}
\caption{Schnakenberg system with periodic boundary conditions at $T=2$: (a) v-RMSE versus time-step size; (b) joint RMSE versus runtime, with each point labeled by its time-step size. All runs use the same $64\times64$ spatial resolution, physical parameters, initial network parameters $\theta_0$.}
\label{fig:tase_dt_comparison}
\end{figure}

To assess the sensitivity of the algorithm to the time-step size, we conduct tests under the same random Stiefel embedding design with $d_{sk}=256$ and $q=8$. As shown in Figure~\ref{fig:tase_dt_comparison}, the accuracy of TASE--TENG improves consistently as the time step is refined. Reducing $\Delta t$ from $5.0\times10^{-2}$ to $2.5\times10^{-2}$ decreases the joint RMSE from $1.9359\times10^{-6}$ to $4.5208\times10^{-7}$; a further reduction to $\Delta t=1.25\times10^{-2}$ lowers the error to $2.5884\times10^{-7}$. The corresponding error reduction factors are approximately $4.28$ and $1.75$, respectively. Hence, although temporal refinement continues to enhance the solution quality, the additional gain becomes noticeably smaller {at the smallest time step, while the runtime increases from $11658\,\mathrm{s}$ to $19028\,\mathrm{s}$.} {This suggests that temporal discretization error is no longer dominant,} with errors from neural manifold projection, optimization, sampling, and spatial discretization becoming increasingly significant.

The effect of the TASE stabilization is more pronounced when the two methods are compared at the same time step $\Delta t=1.25\times10^{-2}$. The joint RMSE of No-TASE--TENG is $7.9896\times10^{-2}$, whereas TASE--TENG achieves $2.5884\times10^{-7}$. Moreover, the joint and demeaned joint RMSEs of TASE--TENG are nearly identical, indicating negligible mean drift. For No-TASE--TENG, however, the error decreases from $7.9896\times10^{-2}$ to $3.1545\times10^{-4}$ after the spatial mean is removed. This large difference shows that a substantial part of the No-TASE--TENG error originates from a global mean drift. These results demonstrate that the TASE stabilization is essential not only for maintaining numerical accuracy but also for suppressing the low-frequency drift generated by the unstabilized target dynamics.
\begin{figure}[!htbp]
\centering
\includegraphics[width=0.96\textwidth]
{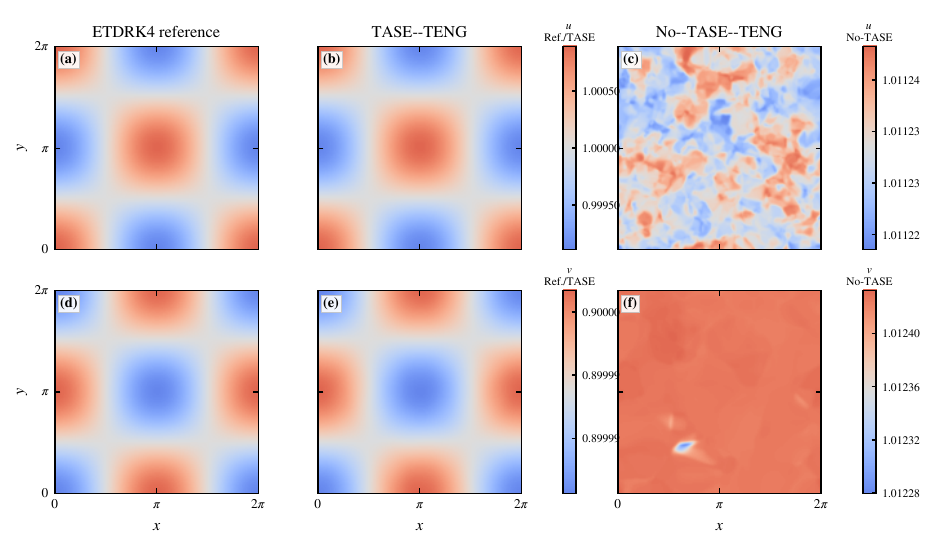}
\caption{
Comparison of the numerical solutions at $T=2.0$ for $\Delta t=1.25\times10^{-2}$. The left, middle, and right columns correspond to the ETDRK4 reference solution, TASE--TENG, and No-TASE--TENG, respectively. }
\label{fig:tase_notase_T2}
\end{figure}
The quantitative comparison in Figure~\ref{fig:tase_dt_comparison} is further supported by the solution profiles shown in Figure~\ref{fig:tase_notase_T2}.
\begin{table}[!htbp]
\centering
\caption{
Resolution-transfer results for the two-dimensional Schnakenberg system with periodic boundary conditions at $T=2$ and $\Delta t=0.05$. The stabilization rule is trained only at $64\times64$ and transferred to the other resolutions without retraining. }
\label{tab:schnakenberg_resolution_transfer}
\setlength{\tabcolsep}{5.0pt}
\renewcommand{\arraystretch}{1.15}
\begin{tabular}{lcccc}
\toprule
Resolution & $u$-RMSE & $v$-RMSE & Joint RMSE & $L^\infty$ \\
\midrule
$48\times48$ & $2.767\times10^{-6}$ & $3.325\times10^{-7}$ & $1.970\times10^{-6}$ & $6.602\times10^{-6}$ \\
$64\times64$ & $2.730 \times10^{-6}$ & $2.030\times10^{-7}$ & $1.935\times10^{-6}$ & $6.772\times10^{-6}$ \\
$80\times80$ & $2.780\times10^{-6}$ & $4.290\times10^{-7}$ & $1.989\times10^{-6}$ & $6.900\times10^{-6}$ \\
\bottomrule
\end{tabular}
\end{table}

At $\Delta t=1.25\times10^{-2}$, the TASE--TENG solution remains in close agreement with the ETDRK4 reference solution at $T=2.0$. In contrast, the No-TASE--TENG solution exhibits a pronounced low-frequency drift, consistent with the substantially larger joint RMSE reported in Figure~\ref{fig:tase_dt_comparison}. The marked reduction from the joint RMSE to the demeaned joint RMSE for No-TASE--TENG further indicates that a significant fraction of this discrepancy is associated with the spatial-mean component of the numerical error.

To examine the resolution transferability of the learned stabilization rule, the rule trained on the $64\times64$ point set is directly {reconstructed} on the $48\times48$ and $80\times80$ point sets without retraining or parameter recalibration. As shown in Table~\ref{tab:schnakenberg_resolution_transfer}, when tested with a time-step size of $\Delta t = 0.05$, the joint RMSE remains around $2\times10^{-6}$ for all three resolutions. Compared with the source $64\times64$ discretization, the joint RMSE increases by only about {$1.8\%$ and $2.8\%$}  after transfer to $48\times48$ and $80\times80$, respectively, while the $L^\infty$ errors stay at essentially the same level. These results indicate that the learned local stabilization rule can be transferred to both coarser and finer spatial resolutions without retraining, while retaining nearly unchanged end-to-end accuracy.

\subsubsection{Homogeneous Neumann boundary conditions}
We next investigate the applicability of the TASE stabilization within the TENG framework under homogeneous Neumann boundary conditions. To isolate the influence of the boundary treatment, all other model parameters and numerical settings are kept unchanged from the periodic case. In particular, the same two-dimensional Schnakenberg system is considered on $\Omega=[0,2\pi]^2$, now subject to the homogeneous Neumann boundary conditions
$$
    \frac{\partial u}{\partial \boldsymbol{n}} = 0, \qquad \frac{\partial v}{\partial \boldsymbol{n}} = 0, \qquad (x,y)\in\partial\Omega,
$$ where $\boldsymbol{n}$ denotes the outward unit normal vector on $\partial\Omega$. Equivalently, on the square domain,
$$
    \begin{aligned} u_x(0,y,t) = u_x(2\pi,y,t) = 0, u_y(x,0,t) = u_y(x,2\pi,t) = 0, \\
     v_x(0,y,t) = v_x(2\pi,y,t) = 0, v_y(x,0,t) = v_y(x,2\pi,t) = 0.\end{aligned}
$$
The initial conditions are also retained from the periodic experiment. For the homogeneous Neumann problem, the reference solution is computed using an  ETDRK4 method combined with a Fourier--cosine spectral discretization. A reference time step of $10^{-4}$ is used. The numerical errors of TASE--TENG and No-TASE--TENG are evaluated on the same $64 \times 64$ spatial points at the final time $T=2$. To isolate the effect of the TASE stabilization, all TENG computations employ the same randomized Stiefel embedding design with $d_{sk}=256$ and $q=8$.

\begin{figure}[!htbp]
\centering
\includegraphics[width=\linewidth]{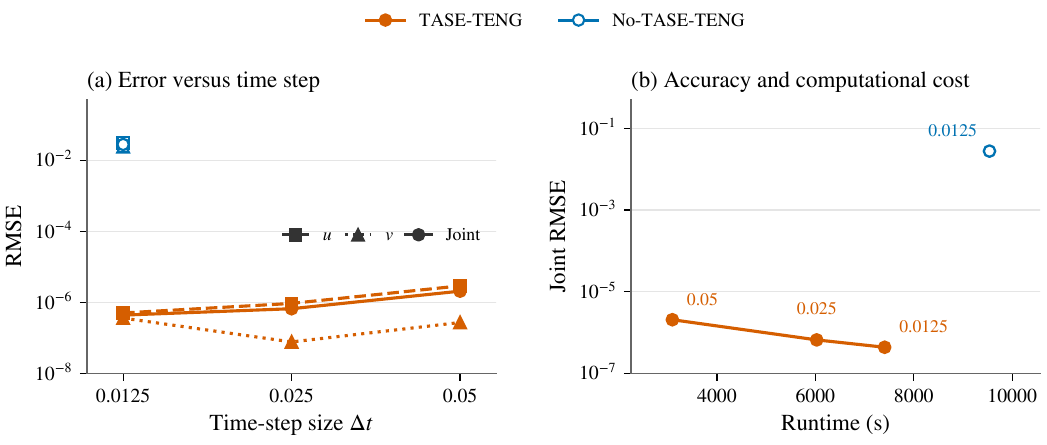}
\caption{Schnakenberg system with homogeneous Neumann boundary conditions at $T=2$: (a) componentwise and joint RMSEs versus time-step size; (b) joint RMSE versus runtime, with each point labeled by its time-step size. The ETDRK4 reference solution uses a Fourier--cosine spectral discretization with $N=256$ and $\Delta t_{\mathrm{ref}}=10^{-4}$. Errors are evaluated on the same $64\times64$ Neumann grid; all TENG runs use $d_{\mathrm{sk}}=256$ and $q=8$.}
\label{fig:schnakenberg_neumann_dt_comparison}
\end{figure}

Figure~\ref{fig:schnakenberg_neumann_dt_comparison} summarizes the accuracy and computational cost of TASE--TENG for three time-step sizes. The joint RMSE decreases monotonically as the time step is refined. Reducing the time step from $\Delta t=5.0\times10^{-2}$ to $2.5\times10^{-2}$ decreases the joint RMSE from $2.0655\times10^{-6}$ to $6.6210\times10^{-7}$, corresponding to a reduction of approximately $67.9\%$. A further reduction to $\Delta t=1.25\times10^{-2}$ decreases the joint RMSE to $4.3714\times10^{-7}$, providing an additional reduction of approximately $34.0\%$. Overall, halving the time step twice reduces the joint RMSE by approximately $78.8\%$.

The componentwise errors do not exhibit identical convergence behavior. The $u$-RMSE decreases consistently from $2.9084\times10^{-6}$ to $9.3319\times10^{-7}$ and then to $5.0249\times10^{-7}$. In contrast, the $v$-RMSE reaches its smallest value, $7.6890\times10^{-8}$, at $\Delta t=2.5\times10^{-2}$ and subsequently increases to $3.6011\times10^{-7}$ at the finest time step. Nevertheless, the reduction in the dominant $u$-component error is sufficiently large that the joint RMSE continues to decrease. This nonuniform componentwise behavior suggests that, {at the smallest time step, the total error is influenced} not only by the temporal discretization but also by the neural-manifold projection, randomized tangent-space approximation, spatial sampling, and nonlinear target-fitting errors. Consequently, the three results should not be interpreted as a clean asymptotic temporal-convergence regime.

The importance of TASE stabilization is evident from the comparison at the common time step $\Delta t=1.25\times10^{-2}$. No-TASE--TENG produces $u$- and $v$-RMSEs of $3.0987\times10^{-2}$ and $2.4096\times10^{-2}$, respectively, resulting in a joint RMSE of $2.7756\times10^{-2}$. At the same time-step size, TASE--TENG achieves a joint RMSE of only $4.3714\times10^{-7}$. Therefore, the No-TASE--TENG error is approximately $6.35\times10^{4}$ times larger. This result demonstrates that randomized parameter-space projection alone is insufficient for this stiff reaction--diffusion problem; stabilizing the Runge--Kutta target before neural-manifold projection is essential for maintaining accuracy.

{The computational results also favor the stabilized method at the smallest tested time step.} TASE--TENG requires approximately $7409\,\mathrm{s}$, whereas No-TASE--TENG requires $9539\,\mathrm{s}$. Thus, the unstabilized computation is approximately $28.7\%$ more expensive. Although applying the TASE operator introduces additional work during target construction, the stabilized computation achieves substantially greater accuracy while also reducing the overall computational time.

\section{Conclusions}
\label{sec:Conclusions}

We have developed TASE--TENG for stiff diffusion-dominated PDEs by incorporating TASE stabilization into the Runge--Kutta target construction before neural-manifold projection. A key component is a learnable graph-based surrogate of the dominant stiff operator. A shared geometry-dependent neural rule generates local interaction coefficients, while a structured factorization enforces discrete self-adjointness, dissipativity, and constant preservation. The surrogate is trained offline and frozen during online integration, where it is used only within the TASE modifier; the PDE right-hand side continues to be evaluated from the continuous neural representation.

For the diffusion and reaction--diffusion problems tested, TASE--TENG reduces numerical errors associated with stiff high-frequency modes and improves the stability and accuracy of neural evolution at the tested time-step sizes. The experiments also demonstrate that the learned local rule can be reused to assemble stabilization operators at different spatial resolutions without retraining. These results support the effectiveness of combining stabilization in physical space with local target realization on the neural manifold.

The present approach uses a time-independent stiff surrogate, whose accuracy may deteriorate when the dominant stiff dynamics vary substantially with the solution or time. Future work will investigate adaptive, structure-preserving update strategies for the surrogate operator, further analyze how operator approximation and neural projection errors affect the stability and accuracy of the coupled method, and extend the proposed framework to the numerical solution of high-dimensional temporally stiff partial differential equations.

\section{Declaration of generative AI and AI-assisted technologies}
{During the preparation of this work, the authors used ChatGPT (OpenAI)(5.6 sol ) to assist with reasoning, selected computational checks, and English language polishing. The authors reviewed all content and take full responsibility for the accuracy and integrity of the entire manuscript.}

\bibliographystyle{alpha}

\begingroup
\setlength{\bibsep}{0pt}
\setlength{\parskip}{0pt}
\bibliography{ref}

@article{r1,
  title={Physics-informed neural networks: A deep learning framework for solving forward and inverse problems involving nonlinear partial differential equations},
  author={Raissi, Maziar and Perdikaris, Paris and Karniadakis, George E},
  journal={ J. Comput. Phys.},
  volume={378},
  pages={686--707},
  year={2019},
  publisher={Elsevier}
}

@article{berman2023randomized,
  title={Randomized sparse neural galerkin schemes for solving evolution equations with deep networks},
  author={Berman, Jules and Peherstorfer, Benjamin},
  journal={Advances in Neural Information Processing Systems},
  volume={36},
  pages={4097--4114},
  year={2023}
}

@article{finzi2023stable,
  title={A stable and scalable method for solving initial value PDEs with neural networks},
  author={Finzi, Marc and Potapczynski, Andres and Choptuik, Matthew and Wilson, Andrew Gordon},
  journal={arXiv preprint arXiv:2304.14994},
  year={2023}
}

@inproceedings{chen2024teng,
  title     = {{TENG}: Time-Evolving Natural Gradient for Solving {PDE}s With Deep Neural Nets Toward Machine Precision},
  author    = {Chen, Zhuo and McCarran, Jacob and Vizcaino, Esteban and Solja{\v{c}}i{\'c}, Marin and Luo, Di},
  booktitle = {Proceedings of the 41st International Conference on Machine Learning},
  series    = {Proceedings of Machine Learning Research},
  volume    = {235},
  pages     = {7143--7162},
  publisher = {PMLR},
  year      = {2024},
  url       = {https://proceedings.mlr.press/v235/chen24ad.html}
}

@article{Ji2021StiffPINN,
  author  = {Ji, Weiqi and Qiu, Weilun and Shi, Zhiyu and Pan, Shaowu and Deng, Sili},
  title   = {Stiff-PINN: Physics-Informed Neural Network for Stiff Chemical Kinetics},
  journal = {The Journal of Physical Chemistry A},
  volume  = {125},
  number  = {36},
  pages   = {8098--8106},
  year    = {2021},
  doi     = {10.1021/acs.jpca.1c05102}
}

@article{Weng2022MPINN,
  author  = {Weng, Yuting and Zhou, Dezhi},
  title   = {Multiscale Physics-Informed Neural Networks for Stiff Chemical Kinetics},
  journal = {The Journal of Physical Chemistry A},
  volume  = {126},
  number  = {45},
  pages   = {8534--8543},
  year    = {2022},
  doi     = {10.1021/acs.jpca.2c06513}
}

@article{Moya2023DAEPINN,
  author  = {Moya, Christian and Lin, Guang},
  title   = {DAE-PINN: A physics-informed neural network model for simulating differential algebraic equations with application to power networks},
  journal = {Neural Computing and Applications},
  volume  = {35},
  pages   = {3789--3804},
  year    = {2023},
  doi     = {10.1007/s00521-022-07886-y}
}

@article{Fabiani2023RanDiffNet,
  author  = {Fabiani, Giacomo and Galaris, Eleftherios and Russo, Lucia and Siettos, Constantinos},
  title   = {Parsimonious physics-informed random projection neural networks for initial value problems of ODEs and index-1 DAEs},
  journal = {Chaos: An Interdisciplinary Journal of Nonlinear Science},
  volume  = {33},
  number  = {4},
  pages   = {043128},
  year    = {2023},
  doi     = {10.1063/5.0135903}
}

@misc{Fabiani2024StabilityRandomProjectionPINN,
  author        = {Fabiani, Giacomo and Bollt, Erik and Siettos, Constantinos and Yannacopoulos, Athanasios N.},
  title         = {Stability of random-projection neural networks in the numerical solution of stiff differential equations},
  year          = {2024},
  eprint        = {2408.15393},
  archivePrefix = {arXiv},
  primaryClass  = {math.NA},
  doi           = {10.48550/arXiv.2408.15393}
}

@article{Zhou2026IPIRNN,
  author  = {Zhou, Hang and Wang, Zhuhong and Qi, Geping and Wang, Yisheng},
  title   = {Integration-based physics-informed randomized neural networks for solving stiff ordinary differential equations},
  journal = {Neurocomputing},
  volume  = {662},
  pages   = {131911},
  year    = {2026},
  doi     = {10.1016/j.neucom.2025.131911}
}

@article{Yao2025SolvingMultiscaleDeePODE,
  title   = {Solving multiscale dynamical systems by deep learning},
  author  = {Yao, Junjie and Yi, Yuxiao and Hang, Liangkai and E, Weinan and Wang, Weizong and Zhang, Yaoyu and Zhang, Tianhan and Xu, Zhi-Qin John},
  journal = {Computer Physics Communications},
  volume  = {316},
  pages   = {109802},
  year    = {2025},
  doi     = {10.1016/j.cpc.2025.109802},
  eprint  = {2401.01220},
  archivePrefix = {arXiv},
  primaryClass  = {cs.LG}
}

@article{Dong2025RandomizedTimeStepping,
  title         = {Randomized time stepping of nonlinearly parametrized solutions of evolution problems},
  author        = {Dong, Yijun and Schwerdtner, Paul and Peherstorfer, Benjamin},
  journal       = {arXiv preprint arXiv:2512.19009},
  year          = {2025},
  eprint        = {2512.19009},
  archivePrefix = {arXiv},
  primaryClass  = {math.NA},
  doi           = {10.48550/arXiv.2512.19009}
}

@article{Conte2026GeneralRKTASE,
  title   = {General Runge--Kutta TASE Methods for Reaction--Diffusion Problems},
  author  = {Conte, Dajana and Montijano, Juan Ignacio and Pagano, Giovanni and Paternoster, Beatrice and R{\'a}ndez, Luis},
  journal = {Journal of Scientific Computing},
  volume  = {106},
  number  = {2},
  pages   = {57},
  year    = {2026},
  doi     = {10.1007/s10915-026-03184-0}
}

@article{kim2025stabilize,
  title={Stabilize physics-informed neural networks for stiff differential equations: Re-spacing layer},
  author={Kim, Eunsuh and Kwon, Heejae and Cho, Sungha and Yeo, Kyongmin and Choi, Minseok},
  journal={Computers \& Mathematics with Applications},
  volume={200},
  pages={167--179},
  year={2025},
  publisher={Elsevier}
}

@article{Eigentler2019Metastability,
  author  = {Eigentler, Lukas and Sherratt, Jonathan A.},
  title   = {Metastability as a Coexistence Mechanism in a Model for Dryland Vegetation Patterns},
  journal = {Bulletin of Mathematical Biology},
  volume  = {81},
  number  = {7},
  pages   = {2290--2322},
  year    = {2019},
  doi     = {10.1007/s11538-019-00606-z}
}

@inproceedings{Gilmer2017MPNN,
  title     = {Neural Message Passing for Quantum Chemistry},
  author    = {Gilmer, Justin and Schoenholz, Samuel S. and Riley, Patrick F. and Vinyals, Oriol and Dahl, George E.},
  booktitle = {Proceedings of the 34th International Conference on Machine Learning},
  pages     = {1263--1272},
  year      = {2017},
  volume    = {70},
  series    = {Proceedings of Machine Learning Research},
  publisher = {PMLR},
  url       = {https://proceedings.mlr.press/v70/gilmer17a.html}
}

@inproceedings{Pfaff2021MeshGraphNets,
  title     = {Learning Mesh-Based Simulation with Graph Networks},
  author    = {Pfaff, Tobias and Fortunato, Meire and Sanchez-Gonzalez, Alvaro and Battaglia, Peter W.},
  booktitle = {International Conference on Learning Representations},
  year      = {2021},
  url       = {https://openreview.net/forum?id=roNqYL0_XP}
}

@misc{raviola2026diracfrenkelonsager,
title        = {A Dirac-Frenkel-Onsager Principle: Instantaneous Residual Minimization with Gauge Momentum for Nonlinear Parametrizations of PDE Solutions},
author       = {Raviola, Matteo and Peherstorfer, Benjamin},
year         = {2026},
eprint       = {2605.00284},
archivePrefix = {arXiv},
primaryClass = {cs.LG},
url          = {https://arxiv.org/abs/2605.00284}
}

@article{wang2021gradient,
  title   = {Understanding and Mitigating Gradient Flow Pathologies
             in Physics-Informed Neural Networks},
  author  = {Wang, Sifan and Teng, Yujun and Perdikaris, Paris},
  journal = {SIAM Journal on Scientific Computing},
  volume  = {43},
  number  = {5},
  pages   = {A3055--A3081},
  year    = {2021},
  doi     = {10.1137/20M1318043}
}

@inproceedings{krishnapriyan2021failure,
  title     = {Characterizing Possible Failure Modes in
               Physics-Informed Neural Networks},
  author    = {Krishnapriyan, Aditi S. and Gholami, Amir and
               Zhe, Shandian and Kirby, Robert M. and
               Mahoney, Michael W.},
  booktitle = {Advances in Neural Information Processing Systems},
  volume    = {34},
  pages     = {26548--26560},
  year      = {2021}
}

@article{bassenne2021tase,
  title   = {Time-Accurate and Highly-Stable Explicit Operators for
             Stiff Differential Equations},
  author  = {Bassenne, Maxime and Fu, Lin and Mani, Ali},
  journal = {Journal of Computational Physics},
  volume  = {424},
  pages   = {109847},
  year    = {2021},
  doi     = {10.1016/j.jcp.2020.109847}
}

@article{calvo2021rktase,
  title   = {A Note on the Stability of Time-Accurate and
             Highly-Stable Explicit Operators for Stiff
             Differential Equations},
  author  = {Calvo, Manuel and Montijano, Jose I. and
             R{\'a}ndez, Luis},
  journal = {Journal of Computational Physics},
  volume  = {436},
  pages   = {110316},
  year    = {2021},
  doi     = {10.1016/j.jcp.2021.110316}
}

@misc{bradbury2018jax,
  title={{JAX}: composable transformations of {Python}+{NumPy} programs},
  author={Bradbury, James and Frostig, Roy and Hawkins, Peter and Johnson, Matthew James and Leary, Chris and Maclaurin, Dougal and Necula, George and Paszke, Adam and VanderPlas, Jake and Wanderman-Milne, Skye and others},
  year={2018}
}

@article{ShiWang2026ConservationTENG,
  author  = {Shi, Zihao and Wang, Dongling},
  title   = {Preserving conservation laws in the time-evolving natural {Galerkin} method via relaxation and projection},
  journal = {Communications in Nonlinear Science and Numerical Simulation},
  volume  = {163},
  pages   = {110743},
  year    = {2026},
  doi     = {10.1016/j.cnsns.2026.110743},
  url     = {https://doi.org/10.1016/j.cnsns.2026.110743}
}

@article{dai2026tinns,
  title={TINNs: Time-Induced Neural Networks for Solving Time-Dependent PDEs},
  author={Dai, Chen-Yang and Chang, Che-Chia and Lin, Te-Sheng and Lai, Ming-Chih and Lai, Chieh-Hsin},
  journal={arXiv preprint arXiv:2601.20361},
  year={2026}
}

@book{Lubich2008,
  author    = {Lubich, Christian},
  title     = {From Quantum to Classical Molecular Dynamics:
               Reduced Models and Numerical Analysis},
  series    = {Zurich Lectures in Advanced Mathematics},
  publisher = {European Mathematical Society},
  address   = {Zurich},
  year      = {2008},
  doi       = {10.4171/067},
  isbn      = {978-3-03719-067-8}
}

@inproceedings{Simonovsky2017ECC,
  author    = {Simonovsky, Martin and Komodakis, Nikos},
  title     = {Dynamic Edge-Conditioned Filters in Convolutional Neural Networks on Graphs},
  booktitle = {2017 IEEE Conference on Computer Vision and Pattern Recognition (CVPR)},
  pages     = {29--38},
  year      = {2017},
  doi       = {10.1109/CVPR.2017.11},
  url       = {https://arxiv.org/abs/1704.02901}
}

@article{Trask2022DDEC,
  author  = {Trask, Nathaniel and Huang, Andy and Hu, Xiaozhe},
  title   = {Enforcing exact physics in scientific machine learning: A data-driven exterior calculus on graphs},
  journal = {Journal of Computational Physics},
  volume  = {456},
  pages   = {110969},
  year    = {2022},
  doi     = {10.1016/j.jcp.2022.110969},
  url     = {https://doi.org/10.1016/j.jcp.2022.110969}
}

@misc{Shaffer2026MEEC,
  author        = {Shaffer, Benjamin D. and Kinch, Brooks and Hsieh, M. Ani and Trask, Nathaniel},
  title         = {A meshfree exterior calculus for generalizable and data-efficient learning of physics from point clouds},
  year          = {2026},
  eprint        = {2605.08436},
  note          = {arXiv preprint arXiv:2605.08436},
  archivePrefix = {arXiv},
  primaryClass  = {cs.LG},
  url           = {https://arxiv.org/abs/2605.08436}
}

@article{eliasof2021pde,
  title={Pde-gcn: Novel architectures for graph neural networks motivated by partial differential equations},
  author={Eliasof, Moshe and Haber, Eldad and Treister, Eran},
  journal={Advances in neural information processing systems},
  volume={34},
  pages={3836--3849},
  year={2021}
}
\endgroup
\appendix
\setcounter{figure}{0}
\section{Target Reachability and Neural-Manifold Compatibility}
{To further investigate the interaction between TASE stabilization and the local neural tangent space, we compare} the No-TASE and TASE targets up to $T=0.5$ with $\Delta t=0.05$. To isolate the effect of the TASE target construction, we additionally evaluate a No-TASE target on the same neural state generated by the TASE trajectory.

\begin{figure}[!htbp]
\centering

\subfloat[Final target-fitting error.\label{fig:app_final_fit}]{\includegraphics[width=0.32\textwidth]{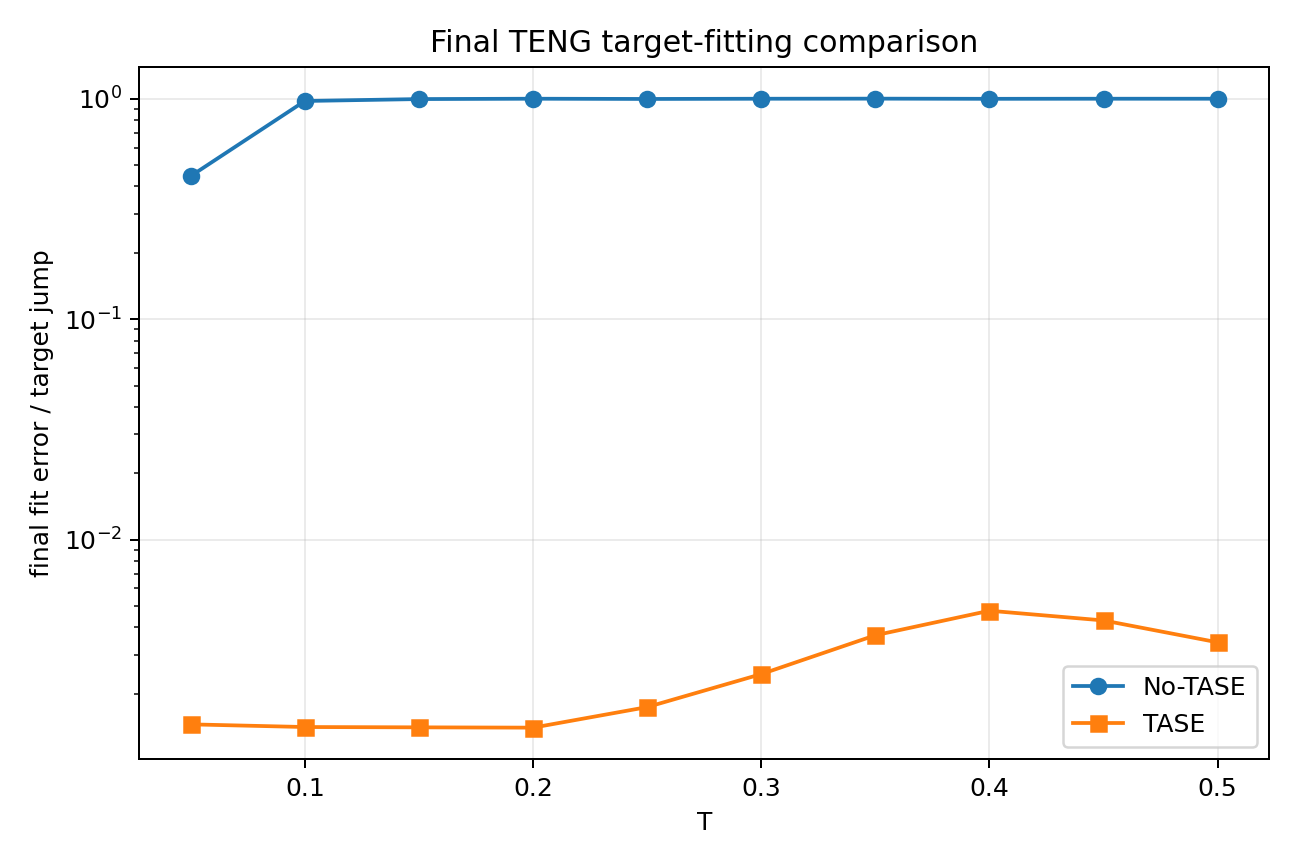}}
\hfill
\subfloat[Stage-2 target reachability.\label{fig:app_stage2_reachability}]{\includegraphics[width=0.32\textwidth]{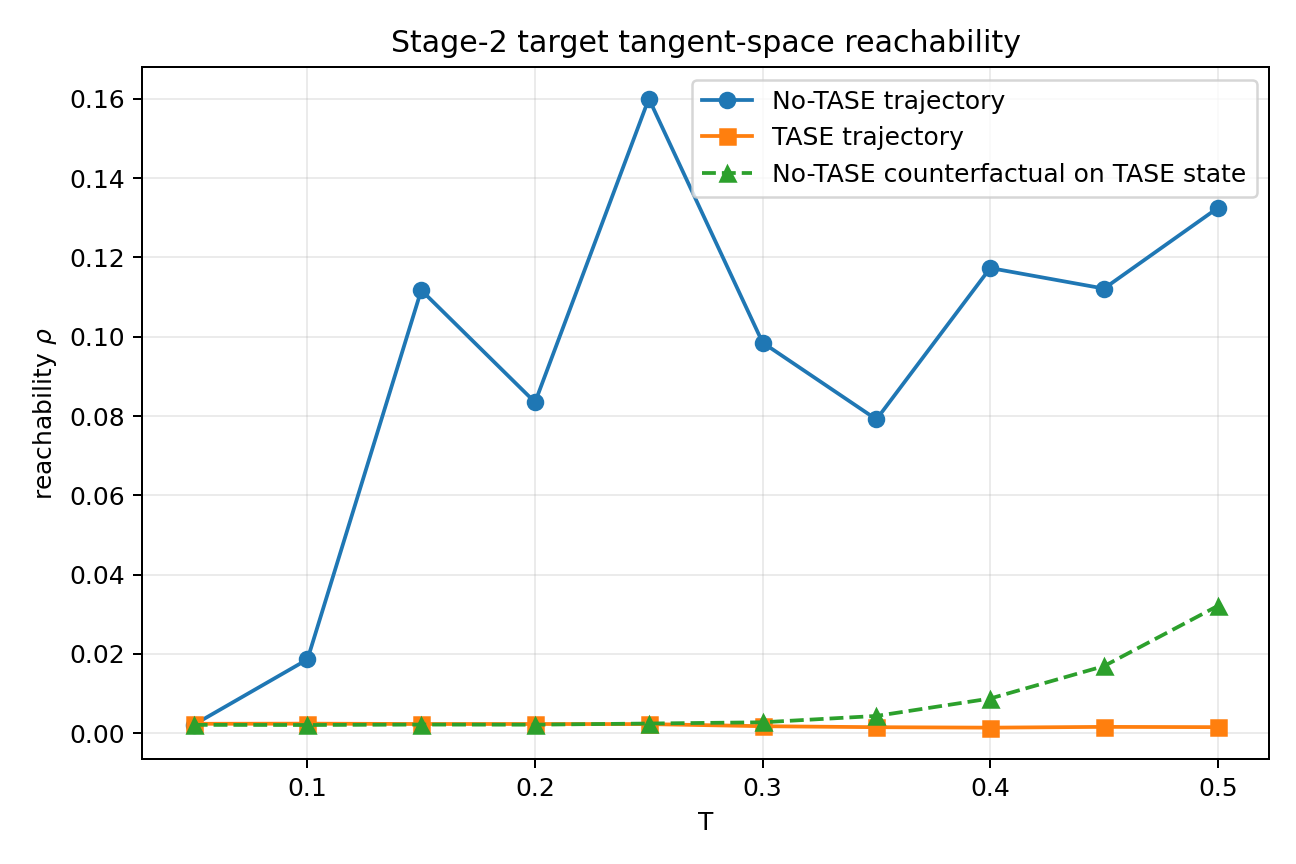}}
\hfill
\subfloat[Final-target reachability.\label{fig:app_final_reachability}]{\includegraphics[width=0.32\textwidth]{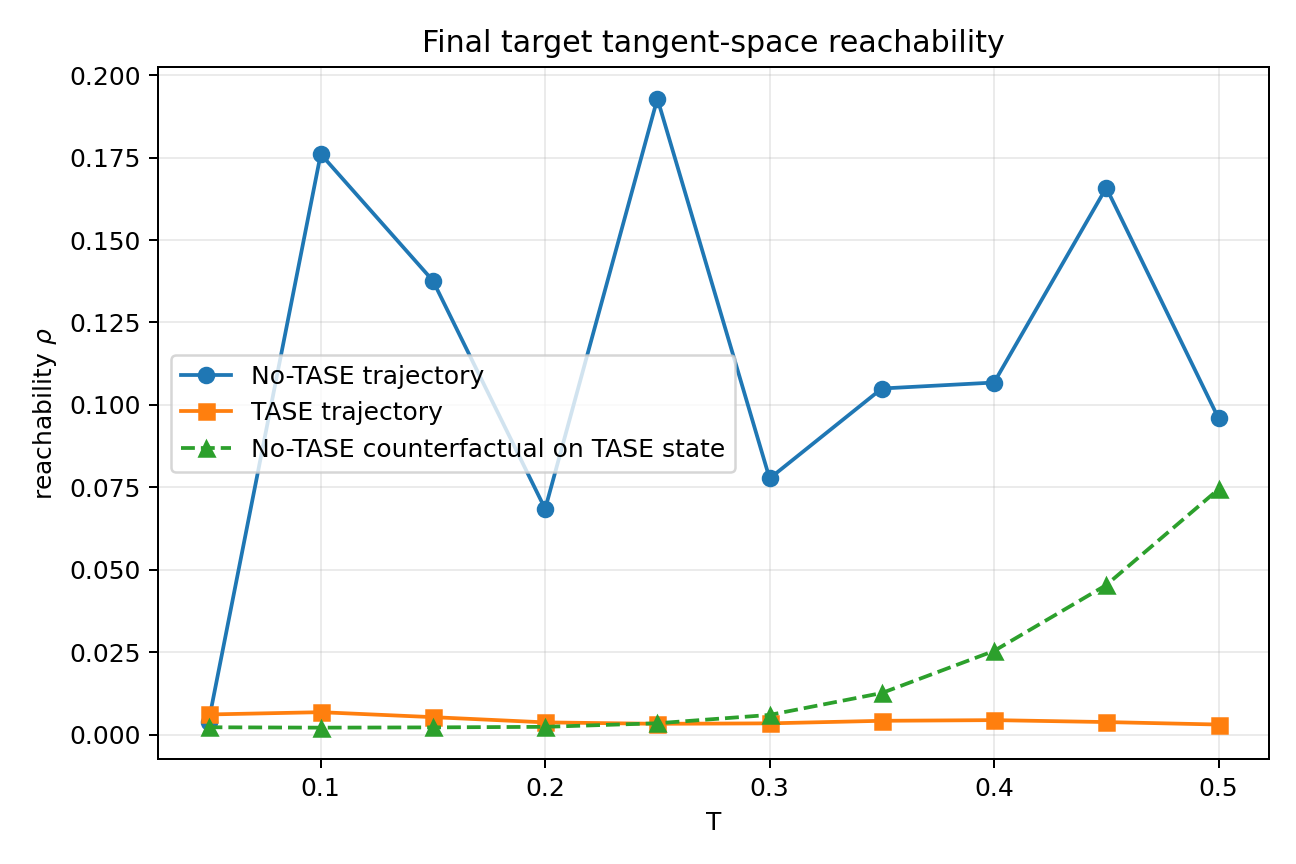}}

\caption{
Comparison of the target-fitting performance and tangent-space reachability between the No-TASE and TASE formulations for the Schnakenberg problem with $\Delta t=0.05$ up to $T=0.5$. The No-TASE target evaluated at the same TASE state is included only as a diagnostic to isolate the effect of the TASE target construction. }
\label{fig:app_tase_reachability}
\end{figure}
Figure~\ref{fig:app_tase_reachability} compares the No-TASE and TASE formulations from three complementary perspectives.
Figure~\ref{fig:app_tase_reachability}(a) reports the relative final target-fitting error,
\[
    \frac{\|u_{\theta_{n+1}}-u_{\mathrm{target}}\|_H}
    {\|u_{\mathrm{target}}-u_{\theta_n}\|_H},
\]
which measures how accurately the TENG iteration represents the time-discrete target.
The No-TASE fitting error rapidly approaches unity, whereas the TASE fitting error remains at the level of $10^{-3}$.
At $T=0.5$, the corresponding values are approximately $9.998\times10^{-1}$ and $3.418\times10^{-3}$.

Figures~\ref{fig:app_tase_reachability}(b) and \ref{fig:app_tase_reachability}(c) report the normalized tangent-space reachability indicator
$$
    \rho = \frac{ \min_{\Delta\theta} \|J_n\Delta\theta-(u_{\mathrm{target}}-u_{\theta_n})\|_H }{ \|u_{\mathrm{target}}-u_{\theta_n}\|_H },
$$
for the Stage-2 target and the final target, respectively. A smaller $\rho$ indicates that the corresponding target increment is better represented by the local neural tangent space.

The blue and orange curves are evaluated along the actual No-TASE and TASE trajectories, respectively. {Since these two trajectories differ after the first time step,} we additionally introduce the green curve as a same-state diagnostic. Specifically, at each state $u_{\theta_n}^{\mathrm{TASE}}$ generated by the TASE trajectory, we keep the neural state and its Jacobian $J_n^{\mathrm{TASE}}$ fixed, but construct the corresponding target with the TASE modification removed. Therefore, the green curve does not represent an additional numerical method; it isolates the effect of the TASE target construction from the accumulated difference between the two trajectories.

The results show that, as the stiff dynamics develop, the TASE target remains substantially more compatible with the local tangent space than the corresponding same-state No-TASE target. At $T=0.5$, for example, the final-target reachability indicators are $3.076\times10^{-3}$ for TASE and $7.442\times10^{-2}$ for the same-state No-TASE target.
\section{Comparison between learned and exact $W_0$}
\label{app:learned_exact_W0}
To further assess the quality of the learned stabilization operator \(W_0^\eta\), we compare it with the exact diffusion operator for the Schnakenberg reaction--diffusion system subject to homogeneous Neumann boundary conditions. The comparison is carried out up to \(T=2\) using a time-step size of \(\Delta t=0.05\). The reference solution is computed using the fourth-order exponential time-differencing Runge--Kutta (ETD--RK4) method with a reference time-step size of \(\Delta t_{\mathrm{ref}}=10^{-4}\). To ensure a fair comparison, all experiments employ the same network architecture, initial network parameters, randomized projection settings, and TENG algorithmic parameters, with the stabilization operator being the only component that is varied.

For the exact case, $W_0$ is constructed from the homogeneous-Neumann
cosine spectral Laplacian,
\begin{equation}
    W_0^{\mathrm{exact}}
    =
    \begin{pmatrix}
        D_u\Delta_N & 0\\
        0 & D_v\Delta_N
    \end{pmatrix},
\end{equation}
\begin{table}[htbp]
\centering
\caption{
Comparison of the learned and exact stabilization operators for the
two-dimensional Schnakenberg system with homogeneous Neumann boundary
conditions at $T=2$ and $\Delta t=0.05, N=64 \times 64 $ .
}
\label{tab:learned_exact_W0}
\begin{tabular}{lccc}
\toprule
$W_0$ & $u$-RMSE & $v$-RMSE & Joint RMSE \\
\midrule
Learned $W_0^\eta$
& $2.9084\times10^{-6}$
& $2.7125\times10^{-7}$
& $2.0655\times10^{-6}$ \\
Exact $W_0^{\mathrm{exact}}$
& $2.8977\times10^{-6}$
& $6.8820\times10^{-7}$
& $2.1059\times10^{-6}$ \\
\bottomrule
\end{tabular}
\end{table}
As shown in Table~\ref{tab:learned_exact_W0}, the learnable operator and the exact operator achieve nearly identical overall computational accuracy. Specifically, the joint RMSE is \(2.0655\times10^{-6}\) with the learnable operator, compared to \(2.1059\times10^{-6}\) with the exact cosine spectral operator. In both cases, the error in the \(u\) component shows almost no difference, while the learnable operator achieves a smaller error in the \(v\) component, which is affected by stiffness.

These results indicate that the learned \(W_0^{\eta^\ast}\) is able to capture the stiffness-related effects required for TASE stabilization with an accuracy comparable to that of the exact diffusion operator. More importantly, this process does not rely on an explicit spectral representation of the underlying differential operator.
\end{document}